\documentclass[11pt]{article}

\usepackage{hyperref}
\usepackage{tikz}
\usepackage{authblk}
\usepackage{enumerate}
\usepackage{mathrsfs}
\usepackage{amsmath}
\usepackage{amssymb}
\usepackage{wrapfig}
\usepackage[thmmarks]{ntheorem}
\usepackage{color}

\renewcommand{\paragraph}{\roman{paragraph}}

\usepackage{hyperref}
\usepackage{tikz}
\usetikzlibrary{arrows,shapes,positioning}
\usetikzlibrary{decorations.markings}
\tikzstyle arrowstyle=[scale=1]
\tikzstyle directed=[postaction={decorate,decoration={markings, mark=at position .65 with {\arrow[arrowstyle]{stealth}}}}]
\tikzstyle reverse directed=[postaction={decorate,decoration={markings, mark=at position .65 with {\arrowreversed[arrowstyle]{stealth};}}}]

\newtheorem{theorem}{Theorem}[section]
\newtheorem{corollary}[theorem]{Corollary}

\newtheorem{conjecture}[theorem]{Conjecture}

\newtheorem{lemma}[theorem]{Lemma}
\newtheorem{proposition}[theorem]{Proposition}
\newtheorem{claim}{Claim}[theorem]
\newtheorem*{LEM2.7}{\textbf{Lemma 2.7}}
\newtheorem*{LEM2.8}{\textbf{Lemma 2.8}}
\newtheorem*{LEM2.9}{\textbf{Lemma 2.9}}
\newtheorem*{LEM2.10}{\textbf{Lemma 2.10}}
\newtheorem*{LEM2.11}{\textbf{Lemma 2.11}}
\newtheorem*{LEM2.12}{\textbf{Lemma 2.12}}
\newtheorem*{LEM2.16}{\textbf{Lemma 2.16}}
\newtheorem*{LEM2.17}{\textbf{Lemma 2.17}}

\newenvironment{proof}{\noindent {\bf
Proof.}}{\rule{3mm}{3mm}\par\medskip}

\newcommand{\D}{\Delta}

\newcommand{\phiv}{\varphi}
\newcommand{\pbar}{\bar{\varphi}}

\newcommand{\CC}{\mathcal{C}}

\newcommand{\JEC}{{\it Europ. J. Combinatorics},  }
\newcommand{\JCTB}{{\it J. Combin. Theory Ser. B.}, }

\newcommand{\JGT}{{\it J. Graph Theory}, }

\newcommand{\GC}{{\it Graph Theory and Combinatorics}, }
\newcommand{\DM}{{\it Discrete Math.}, }

\begin{document}
\title{The average degree of edge chromatic critical graphs with maximum degree seven}
\author{ Yan Cao \\
Scdool of Mathematical Sciences, Dalian University of Technology\\
Dalian, Liaoning, 116024, China\\
Email:  ycao@dlut.edu.cn\\
 Rong Luo\thanks{Partially supported by a grant from  Simons Foundation (No. 839830)} \\
Department of Mathematics,
 West Virginia University\\
 Morgantown, WV 26505 \\
 Email:  rluo@mail.wvu.edu\\
 Zhengke Miao\thanks{Partially supported by NSFC under grant numbers 12031018 and 11971205.  }\\
School of Mathematics and Statistics, Jiangsu Normal University\\
  Xuzhou, Jiangsu, 221116,  China \\
Email: zkmiao@jsnu.edu.cn\\
 Yue Zhao\\ 
Department of Mathematics, University of Central Florida\\
 Orlando, FL 32816-1364\\
 Email: Yue.Zhao@ucf.edu}
 \date{}
\maketitle

\vspace{-1.5cm}

\begin{abstract} In this paper, by developing several new adjacency lemmas about a path on $4$ or $5$ vertices, we show that the average degree of 7-critical graphs is at least 6. It implies Vizing's planar graph conjecture for planar graphs with maximum degree $7$ and its extension to graphs embeddable in a surface with nonnegative Euler characteristic due to Sanders and Zhao (J. Combin. Theory Ser. B 83 (2001) 201-212 and J. Combin. Theory Ser. B 87 (2003) 254-263)
 and Zhang (Graphs and Combinatorics 16 (2000) 467-495).
\end{abstract}

\vspace{-0.2cm}
{\bf Keywords:}. Edge coloring, critical graphs, Euler's formula, planar graphs

\vspace{-0.3cm}
\section{Introduction}
\vspace{-0.3cm}
An {\em edge coloring} of a graph is a function assigning values
(colors) to the edges  of the graph in such a way that any two
adjacent edges receive different colors. A graph is {\em edge
$k$-colorable} if there is an edge coloring of the graph with
colors from $C= \{1,\dots,k\}$. A finite simple graph $G$ of maximum
degree $\Delta$ is {\em class one} if it is edge $\Delta$-colorable.
Otherwise,  $G$ is said to be {\em
class two}, in which case Vizing's Theorem \cite{kn:vizi} guarantees that it is
edge ($\Delta+1$)-colorable.
 $G$ is said to be {\sl edge chromatic critical} (or {\sl critical} for short) if it is connected, class two and $\chi'(G-e) < \chi'(G)$ for every edge $e\in G$. A critical graph $G$  of maximum degree $\Delta$ is called a $\Delta$-{\sl critical graph}.  Vizing  proposed the following conjecture in 1968 \cite{kn:vizing}  on  the average degree of $\Delta$-critical graphs. 

\begin{conjecture}
\label{size conj}
Let $G$ be a $\Delta$-critical graph. Then  $\overline d(G) \geq \Delta - 1 + \frac{3}{|V(G)|}$, where $\overline d(G)$ is the average degree of $G$.
\end{conjecture}

\vspace{-0.3cm}

 There are direct consequences of a progress towards solving this conjecture. For example, if there is a better bound for the size of $\Delta$-critical graphs, then one can obtain better bounds for $\Delta(S)$, where $S$ is a surface and $\Delta(S)=\max\{\Delta(G)|G$ is a class two connected graph that can be embedded in $S\}$. It is well known that if Vizing's conjecture is true for $\Delta=7$, then $\Delta(S) \leq 6$ where $S$ is a surface of  Euler characteristic at least $1$, which was proved in \cite{kn:sz} by other means in 2003. If this average degree conjecture is true, for a  $\Delta$-critical graph $G$, by applying the inequality $\alpha \leq n-\frac{m}{\Delta}$, where $n = |V(G)|$, $m = |E(G)|$, and  $\alpha$ is the independence number of $G$,   one can easily obtain $\alpha \leq \frac{n}{2}$ as $\Delta \rightarrow \infty$. This provides a strong evidence for the independence number conjecture proposed by Vizing in 1968 \cite{kn:vizing},  which claims that if $G$ is a critical graph, then $\alpha \leq \frac{n}{2}$.

Conjecture~\ref{size conj} was verified for $\Delta = 3$ by Jakobsen \cite{kn:jak}, for $\Delta = 4$ by  Fiorini and Wilson \cite{kn:fio},  for $\Delta = 5$ by 
 Kayathri \cite{kn:kaya}, and for $\Delta = 6$ by  Luo, Miao and Zhao \cite{kn:luo}. As for the lower bound of  $\overline d(G)$, 
Woodall \cite{kn:woo} proved that if $G$ is a $\Delta$-critical graph, then $\overline d(G) \geq \frac{2(\Delta +1)}{3}$.    Cao and Chen \cite{Cao2020-1} further improved to $\frac{3\Delta}{4} - 8$ and they \cite{Cao2020-1, Cao2020-2} also showed that Conjecture~\ref{size conj} is asymptotically true.

In this paper, we will prove that if $G$ is a $7$-critical graph, then $\overline d(G) \geq 6$.  This result implies  Vizing's planar graph conjecture for $\Delta = 7$ claiming that every planar graph with maximum degree at least $7$ is class one, which was verified  independently by Sanders and Zhao
\cite{kn:sz} and Zhang \cite{kn:zhanglm} and its extension to graphs embeddable in a surface with nonnegative Euler characteristic due to Sanders and Zhao in  \cite{kn:sz} and \cite{kn:syz}.

Before proceeding, we introduce some notations. Throughout this paper, let $G=(V, E)$ be a simple graph with $n$ vertices, $m$ edges, and maximum degree $\Delta(G)$ (or $\Delta$).  A $k$-{\sl vertex},
$k^+$-{\sl vertex}, or $k^-$-{\sl vertex} is a vertex of degree $k$, at least $k$, or at most $k$, respectively. We use $d(x)$, $d_k(x)$, $d_{k^+}(x)$, $d_{k^-}(x)$ to denote the degree of a vertex $x$, the number of $k$-vertices adjacent to $x$, the number of $k^+$-vertices adjacent to $x$, and the number of $k^-$-vertices adjacent to $x$, respectively. For a vertex $v \in V$, let $N(x) =\{v| xv \in E\}$ be the neighborhood of $v$ in $G$. A \emph{$k$-neighbor}  of a vertex $v$ is a neighbor of $v$ that is a $k$-vertex in $G$, a \emph{$k^+$-neighbor} or \emph{$k^-$-neighbor}  of a vertex $v$ is a neighbor of $v$ that is a $k^+$-{\sl vertex} or $k^-$-{\sl vertex} in $G$. For two disjoint vertex sets $U$ and $U'$, denote by $[U,U']$ the set of edges with one end in $U$ and the other in $U'$.   For a vertex set $A$ of $V(G)$, denote by $N(A) = \cup_{x\in A}N(x)$.

\vspace{-0.3cm}
\section{Lemmas}
\vspace{-0.3cm}
In this section, we present some old  lemmas and develop some  new lemmas needed in the proofs of our main result.

\vspace{-0.3cm}
\subsection{Old lemmas}

\vspace{-0.3cm}

\begin{lemma} (Vizing's Adjacency Lemma \cite{kn:vizi})
\label{Le: VAL}
 Let  $G$ be a $\Delta$-critical graph. Then $d(u) + d(v) \geq \Delta + 2$ for any two adjacent vertices $u$ and $v$, and 
$d_\Delta(x) \geq \max\{2,\Delta-k+1\}$ if $x$ has a $k$-neighbor.
\end{lemma}

 \vspace{-0.5cm}
\begin{lemma} (Luo, Miao, and Zhao \cite{kn:luo})
\label{Le:3-vertex1}
 Let $G$ be a $\Delta$-critical graph with
$\Delta \geq 5$ and $x$ be a 3-vertex.  Then $x$ has  at least two
$\Delta$-neighbors which are not adjacent to any $(\Delta-2)^-$-vertices except $x$.
\end{lemma}

 \vspace{-0.5cm}
\begin{lemma}(Luo, Miao, and Zhao \cite{kn:luo3})
\label{Le:3-vertex2} 
Let $G$ be a $\Delta$-critical graph with
$\Delta \geq 6$ and $x$ be a 3-vertex. Then $x$ has a
$\Delta$-neighbor which is adjacent to at least $\Delta-4 - \lfloor
\frac{\Delta-1}{3}\rfloor$ vertices $z$ with $d(z) = \Delta$ and $d_{(\Delta -3)^-}(z) = 0$.
\end{lemma}

 \vspace{-0.5cm}
\begin{lemma}(Sanders and Zhao
\cite{kn:sz} and Zhang \cite{kn:zhanglm})
\label{Le:delta+2}
 Let $G$ be $\Delta$-critical graph and $xyrs$ be a path with $d(x) + d(y) = \Delta + 2$. Then 
 $d(r) = \Delta$ and $d(s) \geq \Delta -1$. Moreover if $d(x), d(y)<\Delta$, then $d(s) = \Delta$.
 \end{lemma}
 
 \vspace{-0.5cm}
\begin{lemma}(Luo, Miao, and Zhao \cite{kn:luo})
\label{Le:4-vertex} 
Let $G$ be a $\Delta$-critical graph with
$\Delta \geq 6$ and $x$ be a 4-vertex.

\noindent
(1) If $x$ is adjacent to a $(\Delta-2)$-vertex, say $y$,  then  $
N(N(x))\setminus \{x,y\} \subseteq V_{\Delta}$;

\noindent
(2)  If $x$ is not adjacent to any $(\Delta-2)$-vertex and if one of
the neighbors $y$ of $x$ is adjacent to $d(y) - (\Delta -3)$  vertices of degree at most $\Delta-2$, then each of the other three neighbors of $x$
is adjacent to only one $(\Delta -2)^-$-vertex, which is $x$;

\noindent
(3) If $x$ is adjacent to two $(\Delta-1)$-vertices, then each of
the neighbors of $x$ is  adjacent to exactly one $(\Delta
-2)^-$-vertex, which is $x$.
\end{lemma} 

\vspace{-0.3cm}

 The following lemma is a special case of Lemma 2.4  in \cite{kn:sz}  due to Sanders and Zhao. 
 
\vspace{-0.3cm}
 \begin{lemma}
 \label{Cor:375}
 Let $G$ be a $7$-critical graph and $xyz$ be a path in $G$. If $3\leq d(x) \leq 4$, $d(y) = 7$ and $d(x) + d(z) \leq 8$, then $y$ and $z$ have at most $d(x)-3$  common neighbors.
   \end{lemma}
   
\vspace{-0.3cm}
\subsection{New lemmas}
\vspace{-0.3cm}
The following lemmas will be proved in Section 5.
 
 Let $G$ be a $\Delta$-critical graph. For each vertex $v$,  denote
 
 \centerline{$N_{\D\sim 2}(v)=\{z\in N(v):z~has~a~neighbor~of~degree~2\}$}
 
 \vspace{-0.3cm}
\begin{lemma}\label{277}
Let $G$ be a  $\Delta$-critical graph with  $\Delta\ge 7$. Then $|N_{\D\sim 2}(v)|\le 5$ for every $v\in V(G)$.
\end{lemma}

\vspace{-0.5cm}

\begin{lemma}\label{Le:36}
Let $G$ be a $\Delta$-critical graph and $xyrst$ be a path with $d(x)+d(y)=\Delta+2$ and $\max\{d(x),d(y)\}< \Delta$. Then  $d(t)\ge \Delta-2$.
\end{lemma}

\vspace{-0.5cm}
\begin{lemma}
\label{Le:37}
Let $G$ be a $\Delta$-critical graph and $xyrst$ be a path with  $d(x)=3$ and $d(y)=\Delta$.   Suppose that $y$  has a neighbor $z \not \in \{x,r,s\}$ with  $d(z)\le \Delta-2$. Then   $d(s)\ge \Delta-1$; and $d(z)+d(t)\ge \Delta+1$ if $d(t) \le \Delta -4$.
\end{lemma}

So far all adjacency lemmas are about a path on at most four vertices.   Lemma~\ref{Le:37} is  the first  lemma that deals with  a path with five vertices. 

By Lemmas ~\ref{Le:delta+2},  \ref{Le:36}, and \ref{Le:37}, we have the following corollary.

\vspace{-0.3cm}

\begin{corollary}
\label{Cor:2-3-2nd neighbor}
Let $G$ be a $7$-critical graph and $xyrst$ be a path with $d(x) = 3$.  Then we have the following:

\noindent
(1)  if $d(y) = 6$, then $d(r)=d(s) = 7$ and $d(t) \geq 5$.

\noindent
(2) if $d(y) = 7$ and $y$ has another $4^-$-neighbor other than $x$, then  $d(s) \geq 6$ and $d(t) \geq 4$.

\noindent
(3) if $d(y) = 7$ and $y$ has a $5$-neighbor, then  either $d(s) = 6$ and $d(t) \geq 4$ or $d(s) = 7$ and $d(t) \geq 3$.
\end{corollary}

\vspace{-0.5cm}

\begin{lemma}
\label{Le:D+3}
Let $G$ be a $\Delta$-critical graph and  $xy$ be  an edge with $d(x)+d(y)=\Delta+3$ and $\max\{d(x),d(y)\}<\Delta$. Then $x$ has $d(x)-2$ neighbors of degree $\Delta$  having no $(\Delta-2)^-$-neighbors other than $x, y$.
\end{lemma}

\vspace{-0.5cm}
\begin{lemma}
\label{Le:5-vertex}
Let $G$ be a $7$-critical graph and $x$ be a  5-vertex.

\noindent
(1) if $x$ has three 
  6-neighbors, then each $7$-neighbor of $x$ has exactly one $5^-$-neighbor. 
  
  \noindent
  (2)  if $x$  has two 6-neighbors, then   $x$ has  two $7$-neighbors, each of which  has at most two $5^-$-neighbors.
  
  \noindent
  (3) if $x$ has exactly four 7-neighbors,  then $x$ has  two $7$-neighbors, each of which  has at most three $5^-$-neighbors.
\end{lemma}

 \vspace{-0.3cm}

\section{The average degree of 7-critical graphs}

\vspace{-0.3cm}
\subsection{Main result}
In this section we will prove our main result.

\vspace{-0.3cm}
\begin{theorem}
\label{main} 
  $\overline d(G) \geq 6$ for every $7$-critical graph $G$.
\end{theorem}
 \vspace{-0.3cm}
\begin{proof}
Let $G$ be a $7$-critical graph. We define the following subsets of vertices.

$A = \{u| d(u) = 7  \mbox{ and $u$ is adjacent to a $2$-vertex}\}$, 

$B=\{u| d(u) = 6 \mbox{ and $u$ is adjacent to a $3$-vertex}\}$,
 
  $C = \{u| d(u) = 7  \mbox{ and $u$ is adjacent to a $3$-vertex and a $5^-$-vertex}\}$.
  
  The following proposition is straightforward from Lemma~\ref{277} and Corollary~\ref{Cor:2-3-2nd neighbor}.
  
  \vspace{-0.3cm}
   \begin{proposition}
  \label{Prop:ABC}
  Let $x$ be a $7$-vertex  which is not adjacent to a $5^-$-vertex. Then at most one of the three sets $N(x)\cap A$, $N(x)\cap B$, and $N(x)\cap C$ is a nonempty set. Moreover $|N(x)\cap A|\leq 5$ and $|N(x)\cap B|\leq 1$.
  \end{proposition}

\vspace{-0.3cm}
For each vertex $x$, denote by $M(x) = d(x) -6$ to be the initial charge of $x$.

\noindent
{\bf R1}    Let $u$ be a $7$-vertex  not adjacent to a $5^-$-vertex but adjacent to a vertex in $A\cup B \cup C$. Then $u$ sends $\frac{1}{|N(x)\cap A| + |N(x)\cap B|+|N(x)\cap C|}$ to each neighbor   in $A\cup B \cup C$.

\medskip \noindent
{\bf R2} Let $u$ be  a $7$-vertex  adjacent to a $5^-$-vertex. Then $u$ sends $\frac{1}{d_{5^-}(u)}$ to each neighbor with degree $4$ or $5$, $1$ to each $3$-neighbor, and $2$ to each $2$-neighbor.

\medskip \noindent
{\bf R3} Every $6$-vertex sends $1$ to each $3$-neighbor.

\medskip \noindent
{\bf R4} If a $5$-vertex $u$ is adjacent to a $7$-vertex $v \in C$, then $u$ sends $\frac{1}{8}$ to  $v$.

\medskip \noindent
{\bf R5} If a $4$-vertex is adjacent to a $5$-vertex, then the $4$-vertex receives $\frac{1}{2}$ from its $5$-neighbor.

\medskip

Denote by $M'(x)$ to be the new charge of the vertex $x$.  We have the following estimation for $M'(x)$. 

\begin{enumerate}
\item[(I)] 
\label{I} Let $u$ be a vertex with degree $2$ or $3$. Then $M'(u) =0$. 

By (R2), each $2$-vertex receives $2$ from each neighbor. By Lemma~\ref{Le: VAL},  each $3$-vertex is not adjacent to a $5^-$-vertex. Thus by (R2) and (R3), each $3$-vertex receives $1$ from each neighbor. Therefore $M'(u) = 0$ if $d(u) = 2$ or $3$.

\item[(II)] 
\label{II}
 Let $uv$ be an edge with $d(u) + d(v) = \Delta + 2 = 9$ and $3\leq d(u) \leq d(v) < 7$.
  Then  $M'(u) \geq 0$ and $M'(v) \geq 1$.

Let $w\in N(u)\cup N(v)$ and $w\not \in \{u,v\}$. If  $w\in N(u)\cap N(v)$,  then by  Lemma~\ref{Le:delta+2},  $d(w) = 7$,  and  $w$  has only two  $6^-$-neighbors. Thus by (R2), $w$ sends $\frac{1}{2}$ to each of $u$ and $v$ if $d(u) = 4$ and $d(v) = 5$ and sends $1$ to $u$, $0$ to $v$ if $d(u) = 3$ and $d(v) = 6$.

If $w \not \in N(v) \cap N(u)$, then by  Lemma~\ref{Le:delta+2},  $d(w) = 7$ and  $w$  has only one $6^-$-neighbor, which is  either $u$ or $v$. If $w \in N(u)$, then by (R2), $w$ sends $1$ to $u$.  Assume  $w\in N(v)$.  If $d(v) = 6$, then     $v\in B$, and by Proposition~\ref{Prop:ABC}, $w$ sends $1$ to $v$. If $d(v) = 5$, then $N(w)\cap (A\cup B\cup C)=\emptyset$ by Lemma~\ref{Le:36} and  thus $w$ sends $1$ to $v$ by (R2).   Therefore in any case $w$ sends $1$ to either $u$ or $v$ if $w \not \in N(v) \cap N(u)$.

If $d(u) = 4$ and $d(v) = 5$, then $u$ receives $\frac{1}{2}$ from $v$ by (R2). Thus $M'(u) \geq 4-6 + 4\times \frac{1}{2}= 0$ and $M'(v)=5-6+\frac{1}{2}|N(u)\cap N(v)|+|N(v)\setminus N(u)|-\frac{1}{2}\ge 5-6+\frac{3}{2}+1-\frac{1}{2}\ge 1$.

 If $d(u) = 3$ and $d(v) = 6$, then  $M'(u) = 0$ by (I) and  $v$ sends  $1$ to $u$ by (R3). Thus $M'(v)=6-6+|N(v)\setminus N(u)|-1\ge 6-6+3-1> 1$.

\item[(III)] 
\label{III}
 Let $u$ be a $4$-vertex with four $6^+$-neighbors. Then $M'(u) > 0 $ unless $u$ has either four $7$-neighbors or has two $6$-neighbors and two $7$-neighbors, in which case $M'(u) \geq 0$.

By Lemma~\ref{Le: VAL},  $u$ is adjacent to at least two $7$-vertices and each $7$-neighbor of  $u$ is adjacent to at most three  $5^-$-vertices. 

If $u$ has a $7$-neighbor $v$ adjacent to  three $5^-$-vertices, then by Lemma~\ref{Le:4-vertex},  $u$ is adjacent to four $7$-vertices and except $v$, each $7$-neighbor is adjacent to only one $5^-$-vertex. Therefore by (R2),  $M'(u) \geq 4-6 + 3\times 1 + \frac{1}{3} = \frac{4}{3}$.

Now assume that each $7$-neighbor is adjacent to at most two $5^-$-vertices. Then $u$ receives at least $\frac{1}{2}$ from each $7$-neighbor.

If $u$ has four $7$-neighbors, then $M'(u) \geq 4-6 + 4\times \frac{1}{2} = 0$.

If $u$ has a $6$-neighbor, then by Lemma \ref{Le:D+3}, there are two  $7$-neighbors of $u$ having only one $5^-$-neighbor. Thus $M'(u) \geq -2 + 2 + \frac{1}{2}(d_7(u) -2) \geq 0$ with equality when $u$ has exactly two $6$-neighbors and two $7$-neighbors.

\item[(IV)] 
\label{IV}  $M'(u) > 0$ for  each $5$-vertex $u$ with five $5^+$-neighbors.

By Lemma~\ref{Le: VAL},  $u$ is adjacent to at least two $7$-vertices and each $7$-neighbor of  $u$ is adjacent to at most four  $6^-$-vertices. 

If $v$ is a $7$-neighbor of $u$ and $v$ is adjacent to a $3$-vertex, then $v$ sends $\frac{1}{2}$ to $u$ by (R2) and $u$ sends $\frac{1}{8}$ to $v$ by (R4). Therefore the total net charge $u$ receives from $v$ is $\frac{3}{8}$. 

Thus in general, $u$ receives at least $\min\{\frac{3}{8},\frac{1}{4}\}$ from each $7$-neighbor. 

If $u$ has  at least four  $7$-neighbors, then  by Lemma~\ref{Le:5-vertex}(3),  $M'(u)\geq -1+ 2\times\frac{1}{4}+2\times \frac{1}{3}>0$.

Now assume that $u$ is adjacent to at most three $7$-vertices. 

If $u$ is adjacent to a $5$-vertex, then by  Lemma~\ref{Le:D+3}, $u$ has  three $7$-neighbors, each of which  could be adjacent to at most two $5^-$-vertex ($u$ and the $5$-neighbor of $u$). Thus $M'(u) \geq -1 + 3\times \frac{1}{2} = \frac{1}{2} > 0$.

Finally, we may assume that $u$ is adjacent to at least two $6$-vertices and at most three $7$-vertices. By Lemma~\ref{Le:5-vertex}(1) and (2),  $M'(u) \geq -1 +  \min\{\frac{1}{4} + 2\times \frac{1}{2}, 1 + 1\} > 0$.

\item[(V)] 
\label{V} 
Let $u$ be a $6$-vertex  adjacent to six $4^+$-vertices. Then by the discharging rules, $M'(u) = M(u) = 0$.

\item[(VI)] 
\label{VI}
 $M'(u) \geq 0$ if $d(u) = 7$.
 
 Let  $u$ be a $7$-vertex.  Then $u \not \in B$. By (R1) and (R2),  we have $M'(u) \geq  0$ if $u  \not \in A\cup  C$.

(a) Assume $u \in A$ (that is $u$   has a $2$-neighbor $v$). 

Let $w$ be the other neighbor of $v$ and $x \in N(u)\setminus \{v,w\}$. Then by Lemma~\ref{Le:delta+2}, $d(x) = 7$ and  $x$ is not adjacent to a $5^-$-vertex.  Since $u\in A$, by Proposition~\ref{Prop:ABC}, $x$ is adjacent to at most five vertices in $A\cup C$. Thus by (R1), $x$ sends at least $\frac{1}{5}$ to $u$. Since  $|N(u)\setminus \{v,w\}| \geq 5$, we have $M'(u) \geq 7-6-2 + 5\times \frac{1}{5}=0$.

(b) Assume  $u \in C$ (that is $u$  is adjacent to a $3$-vertex  $x$ and another $5^-$-vertex $z$).

By Lemma~\ref{Le: VAL}, $x$ and $z$ are not adjacent and $u$ has five $7$-neighbors. By Lemma~\ref{Cor:375}, $u$ and $z$ have no common neighbor. Thus  $u$ has at least three $7$-neighbors which are not adjacent to $x$ or $z$.  Let $w$ be such a $7$-neighbor of $u$. By Proposition~\ref{Prop:ABC}, $N(w) \cap (A\cup B)= \emptyset$.

If $d(z) \leq 4$, then $3 \leq d(z) \leq 4$ by Lemma~\ref{Le: VAL}, and thus $u$ sends at most $1$ to each of $x$ and $z$. 
By Corollary~\ref{Cor:2-3-2nd neighbor}(2), $u$ is the only vertex in $C$ adjacent to $w$. So $w$ sends $1$ to $u$ by (R1).  Thus $M'(u) \geq 7-6 -1-1 + 3 =1$.

If $d(z) = 5$, then $w$ is adjacent to at most seven vertices in $C$ and thus sends at least $\frac{1}{7}$ to $u$ by (R1). By (R2), $u$ sends $1$ to $x$ and $\frac{1}{2}$ to $z$ and by (R4), $z$ sends $\frac{1}{8}$ to $u$. Therefore $M'(u) \geq 7-6 -1-\frac{1}{2} + \frac{1}{8} + \frac{3}{7} > 0$. This completes the proof of (VI).
\end{enumerate}

By (I)-(VI), $M'(x) \geq 0$ for each vertex $x$ and thus $0 \leq \sum_{x\in V} M'(x) = \sum_{x\in V} M(x) = (\overline d(G)-6)|V|$. Therefore $\overline d(G) \geq 6$. This completes the proof of the theorem.
\end{proof}

\vspace{-0.3cm}
\subsection{Concluding remarks}
\vspace{-0.3cm}
\noindent
One may wonder why our result does not include the term $\frac{3}{|V|}$ in the lower bound for the average degree as Conjecture~\ref{size conj} states. The reason is that we can construct some infinite families of graphs with maximum degree $7$ and average degree $6$  which satisfy all currently known adjacency lemmas.  For example, for any positive integer $t$, consider a graph $G$  with degree sequence $(4^t, 7^{2t})$ such that each $4$-vertex is adjacent to four $7$-vertices and each $7$-vertex is adjacent two $4$-vertices. One can easily check that $G$ satisfies all adjacent lemmas that we currently have and  $\overline d(G) = 7-1=6$.  The above example can be generalized  for arbitrary maximum degree $\Delta = 2k+1 \geq 7$. For each $t \geq 1$, let $G$ be a graph with degree sequence $(k^t, \Delta^{kt})$ such that each $k$-vertex is adjacent to $k$ vertices of degree $\Delta$ and each $\Delta$-vertex is adjacent to  exactly one $k$-vertex. Then $\overline d(G) = \Delta-1=2k$ and  $G$ satisfies all adjacency lemmas that we know.  

The above examples and several other examples not only present a challenge but also indicate the necessity to develop new adjacency lemmas to attack Conjecture~\ref{size conj} and other edge coloring problems. In particular, so far all adjacency lemmas are about a path on at most four vertices.   Lemma~\ref{Le:37} is indeed a lemma that deals with  a path with five vertices and it is the key lemma in the proof of our main result, but it is only for degree $3$-vertices. To completely solve the case of $7$-critical graphs and beyond,   more general adjacency lemmas concerning paths on five  vertices are needed although it is very challenging to develop such lemmas.   It  would be  practical and very useful  to use computer program to  complete  the remaining cases for 7-critical graphs and to develop some forbidden structures for critical graphs in general.

 \vspace{-0.3cm}

\section{Applications to graphs embedded on surfaces with nonnegative Euler characteristics}
 \vspace{-0.3cm}
Theorem~\ref{main}  clearly implies that every planar graph with maximum degree  $7$ is class one  which was conjectured by Vizing and  independently proved by Sanders and Zhao \cite{kn:sz}, and Zhang~\cite{kn:zhanglm} and its extension to projective planar graphs \cite{kn:syz} since every graph  which can be embedded in a  plane or a projective plane has average degree strictly less than $6$.  Our result also implies the following result due to Sanders and Zhao~\cite{kn:syz}. 

\vspace{-0.3cm}
\begin{theorem} (Sanders and Zhao~\cite{kn:syz})
Let $G$ be a graph with maximum degree $7$. If $G$ can be embedded in the torus or Klein bottle, then $G$ is class one.
\end{theorem}

 \vspace{-0.3cm}
\begin{proof} Prove by contradiction. Suppose that $G$ is not class one. Then we may assume that  $G$ is $7$-critical.  By Euler's formula,  $\overline d(G) \leq 6$. By Theorem~\ref{main},  we have  $\overline d(G) = 6$ and $d(f) = 3$ for each face $f$.  Since $G$ is simple, we further have $\delta \geq 3$. Denote by $M'(x)$ the new charge of the vertex $x$ and $A,B,C$ the sets defined  in  the previous section.  Then $\sum_{x\in V(G)}M'(x) = \sum_{x\in V(G)} (\overline d(G)-6) = 0$. Thus  $M'(x) = 0$ for every vertex $x$ in $G$.   
 
 Since $\delta (G) \geq 3$, we have $A = \emptyset$. By (II) and (IV) in the proof of Theorem~\ref{main}, $d(u) + d(v) \geq \Delta + 3$ and there are no $5$-vertices in $G$.   Thus $B = \emptyset$.  Since every face is a $3$-face and $G$ is $2$-connected, every two adjacent vertices share at least two common neighbors. 
 \begin{claim}
 \label{CL:3,4}
  $\delta (G) = 4$ and every $4$-vertex is adjacent to exactly two $7$-vertices and two $6$-vertices.
   \end{claim}
    \vspace{-0.3cm}
 \begin{proof}
  Let $y$ be a $7$-vertex  with a neighbor $x$ where $3\leq d(x) \leq 4$.  Since any two adjacent vertices share at least two neighbors,  by Lemma~\ref{Cor:375}, $y$ is adjacent to only one $4^-$-vertex. Since there are no $5$-vertices in $G$, $y$ is adjacent to exactly one $5^-$-vertex.  
   This implies $C = \emptyset$. Therefore $A= B = C = \emptyset$. Hence every $7$-vertex is adjacent to a $4^-$-vertex otherwise $M'(x) = M(x) = 1 > 0$ if $x$ is a $7$-vertex without a $4^-$-neighbor. Therefore every $7$-vertex has exactly one $4^-$-neighbor.
  
  If there is a $3$-vertex, by Lemma~\ref{Le:3-vertex2}, there is one $7$-vertex $x$ that has no $4^-$-neighbors,  a contradiction.  Therefore $\delta = 4$ and every $7$-vertex is adjacent to exactly one $4$-vertex. By (III), every $4$-vertex is adjacent to exactly two $7$-vertices and two $6$-vertices.
 \end{proof}
 
 Denote by $V_i$ the set of $i$-vertices and $n_i = |V_i|$. Then by Claim~\ref{CL:3,4}, $n_4 = 2n_7$ and $n_4\leq 2n_6$.
 
 Since every $7$-vertex is adjacent to a $4$-vertex, every $7$-vertex is adjacent to at least $4$ vertices in $V_7$ and every vertex has at least two neighbors in $V_7$ by Lemma~\ref{Le: VAL}. Thus $ 2n_6 + 2n_4 \leq |[V_7, V_6\cup V_4]| \leq 3n_7$. This implies 
 $6n_7 = 3n_4 \leq 3n_7$. This contradiction completes the proof of the theorem.
 \end{proof}

 \vspace{-0.5cm}

\section{Proofs of  new lemmas}
\label{SEC: proof-lemma}
 \vspace{-0.3cm}

Before giving the proofs, we first introduce some notations and  lemmas that are needed in this section.

The set of all $k$-edge-colorings of a graph $G$
is denoted by $\mathcal{C}^k(G)$. Let $\varphi\in \mathcal{C}^k(G)$.  For any color $\alpha$, let
$E_{\alpha} =\{e\in E : \phiv(e) =\alpha\}$. For any two distinct colors $\alpha$ and $\beta$, denote
by $G_{\varphi}(\alpha,\beta)$
 the subgraph of $G$ induced by $E_{\alpha} \cup E_{\beta}$. The components of $G_{\varphi}(\alpha,\beta)$ are called
{\it $(\alpha,\beta)$-chains}.  Clearly, each $(\alpha,\beta)$-chain is either a path or a cycle of edges  alternately colored with $\alpha$ and $\beta$.
For each $(\alpha, \beta)$-chain $P$, let $\varphi/P$ denote the $k$-edge-coloring obtained from $\varphi$ by exchanging colors $\alpha$ and $\beta$ on $P$.

For any $v\in V$, let $P_v(\alpha, \beta, \varphi)$ denote
the unique $(\alpha, \beta)$-chain containing $v$. Notice that, for any two vertices $u, \, v\in V$, either $P_u(\alpha, \beta, \varphi)=P_v(\alpha, \beta, \varphi)$
or $P_u(\alpha, \beta, \varphi)$ is vertex-disjoint from $P_v(\alpha, \beta, \varphi)$. This fact will be used very often without mentioning. For convenience, we define $P_v(\alpha,\beta,\phiv)=v$ and $\phiv/P_v(\alpha,\beta,\phiv)=\phiv$ when $\alpha=\beta$.

For any $v\in V$, let $\varphi(v)=\{\varphi(e)\,: e\in E(v)\}$ denote the set of colors presented at $v$  and
$\pbar(v) = C\setminus \varphi(v)$ the set of colors not assigned to any edge incident to $v$, which are called  {\it missing} colors at $v$. For a vertex set $X\subseteq V(G)$, we call $X$  {\it elementary} (with respect to $\varphi$) if all missing color sets $\pbar(x)$ ($x\in X$) are mutually  disjoint.

A \emph{multi-fan} at $x$ with respect to the edge $e=xy \in E(G)$ and the coloring $\phiv\in \mathcal{C}^\D(G-e)$ is a sequence $F = (x, e_1,y_1,\ldots, e_p,y_p)$ with $p \ge 1$ consisting of edges $e_1,e_2,\ldots,e_p$ and vertices $x,y_1,y_2,\ldots,y_p$ satisfying the following two conditions:
\vspace{0.1 cm}

$\bullet$ The edges $e_1,e_2, \ldots, e_p$ are distinct, $e_1=e$ and $e_i =  xy_i$ for $i=1, \ldots, p$.

\vspace{0.1 cm}

$\bullet$ For every edge $e_i$ with $2\le i\le p$, there is a vertex $y_j$ with $1\le j<i$ such that $\varphi(e_i)\in \pbar(y_j)$.

Note that a multi-fan is slightly more general than a Vizing-fan which requires $j=i-1$ in the second condition.

\begin{lemma}{(Stiebitz, Scheide, Toft and Favrholdt~\cite{tec})}
\label{vf}
 Let $G$ be a $\Delta$-critical graph, $xy_1=e\in E(G)$ and $\phiv\in \mathcal{C}^\D(G-e)$. If $F = (x, e_1, y_1, \ldots, e_p, y_p)$ is a multi-fan at $x$ with respect to $e$ and $\varphi$. Then the following statements hold:

(a) $\{x,y_1,y_2,\ldots,y_p\}$ is elementary.

(b) If $\alpha \in \pbar(x)$ and $\beta \in \pbar(y_i)$ for some $i$, then $P_x(\alpha,\beta,\phiv)=P_{y_i}(\alpha,\beta,\phiv)$.
\end{lemma}

\vspace{-0.3 cm}

 The following  lemma is a direct corollary of  Lemma~\ref{vf}.

 \begin{lemma}
 \label{path}
 Let $G$ be a $\Delta$-critical graph, $xy=e\in E(G)$ and $\phiv\in \mathcal{C}^\D(G-e)$.  Let $xyz$ be a path.
 
 (1)  If $d(z) \leq 2\Delta - (d(x)+d(y)) + 1$, then $\alpha =\phiv(yz) \in \phiv(x)\cap \phiv(y)$ and for any color $\beta \in \pbar(z)\cap (\pbar(x)\cup \pbar(y))$, $P_z(\alpha, \beta, \phiv)$ ends at $x$ or $y$.
 
 (2) If $\phiv(yz) \in \pbar(x)$, then $\pbar(x) \cup \pbar(y) \subseteq \phiv(z)$ and thus $d(z) \geq 2\Delta - (d(x)+d(y)) + 2$.
 \end{lemma}

\vspace{-0.3 cm}
A \emph{Kierstead path} with respect to $e=y_0y_1$ and $\phiv\in \CC^\D(G-e)$ is a path $K = y_0y_1\cdots y_p$ with $p \ge 1$ such that for every edge $y_iy_{i+1}$ with $1\le i\le p-1$, there is a vertex $y_j$ with $0\le j<i$ such that $\varphi(y_iy_{i+1})\in \pbar(y_j)$.

Clearly a Kierstead path with $3$ vertices is a multi-fan with center $y_1$.  The next two lemmas are elementary properties of  a Kierstead path with $4$ vertices.

\begin{lemma}\label{p4}{(Kostochka and Stiebitz~\cite{tec}, Luo and Zhao~\cite{LuoZhao})}
Let $G$ be a $\Delta$-critical graph, $y_0y_1=e\in E(G)$ and $\phiv\in \mathcal{C}^\D(G-e)$.  Let $K=y_0y_1y_2y_3$ be a Kierstead path with respect to $e$ and $\varphi$.
Then $V(K)$ is elementary unless $d(y_1)=d(y_2)=\D(G)$, in which case, all colors in $\pbar(y_0),\pbar(y_1),\pbar(y_2)$ and $\pbar(y_3)$ are distinct except one possible common missing color in $\pbar(y_3)\cap(\pbar(y_0)\cup\pbar(y_1))$.  
\end{lemma}

\vspace{-0.5 cm}

\begin{lemma}
\label{p4link}
Let $G$ be a $\Delta$-critical graph, $y_0y_1=e\in E(G)$ and $\phiv\in \mathcal{C}^\D(G-e)$. Suppose that $K=y_0y_1y_2y_3$ is a Kierstead path with respect to $e$ and $\varphi$, $\min\{d(y_1),d(y_2)\}<\D$, $\alpha\in \pbar(y_3)$ and $\beta\in \pbar(y_i)$  for some $i \in \{0,1,2\}$. If $\beta\notin \{\phiv(y_1y_2),\phiv(y_2y_3)\}$, then $P_{y_3}(\alpha, \beta, \phiv)$ ends at $y_i$.
\end{lemma}
\vspace{-0.3 cm}
\begin{proof}
Since $K$ is a Kierstead path and $\{y_0,y_1,y_2,y_3\}$ is elementary by Lemma~\ref{p4}, we have $\alpha\notin \{\phiv(y_1y_2),\phiv(y_2y_3)\}$. Suppose to the contrary that $P_{y_3}(\alpha, \beta, \phiv)$ does not end at $y_i$. Then after interchanging $\alpha,\beta$ on this path, $K$ is still a Kierstead path, but $\beta$ is missing at both $y_i$ and $y_3$, a contradiction to Lemma~\ref{p4}. This completes the proof.
\end{proof}

\begin{figure}[!ht]
    \centering
    \includegraphics[width=5.1in,height=3.4in]{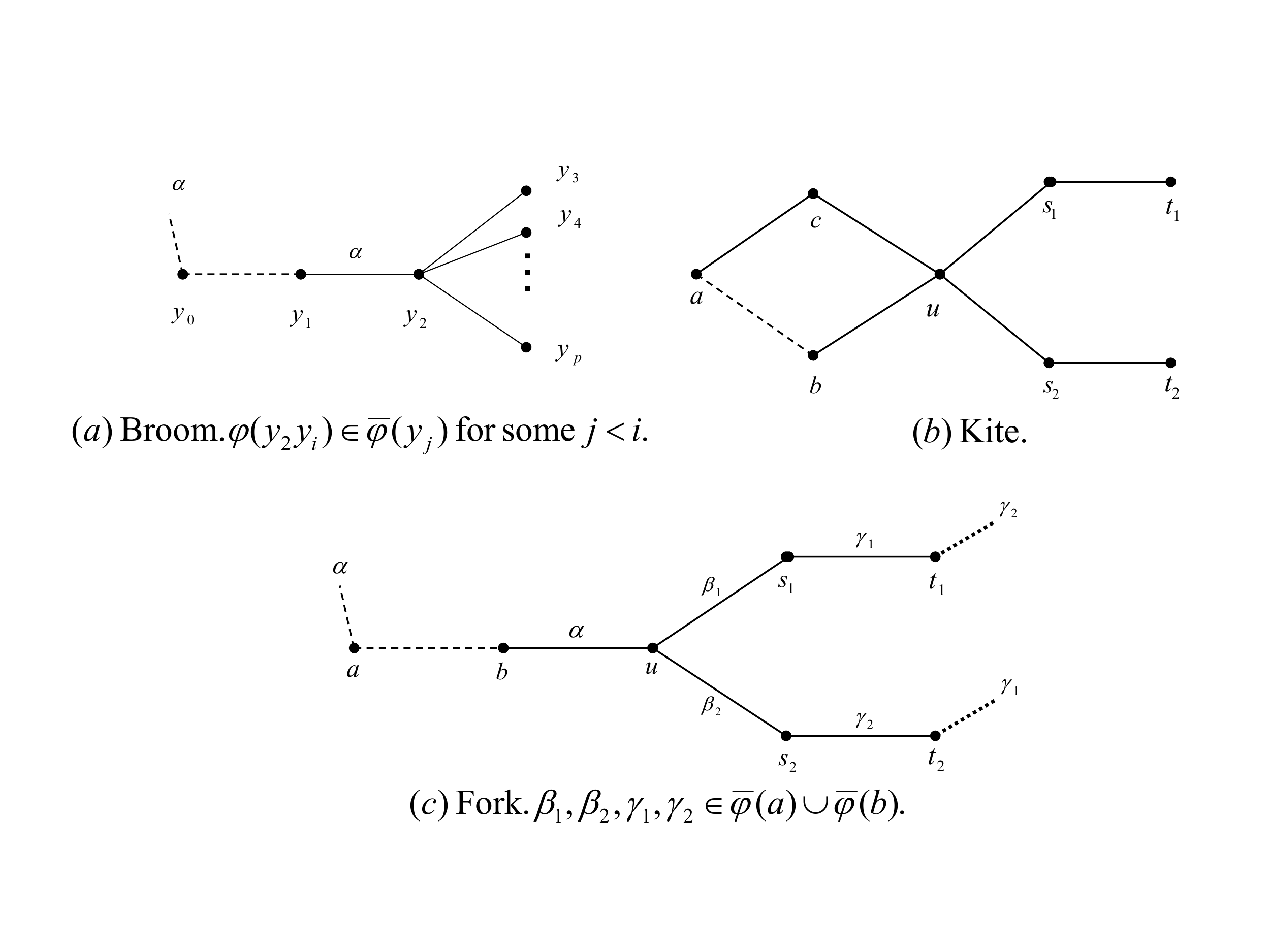}
	\caption{Brooms, kites and forks.}\label{fig}
\end{figure}

A \emph{$\phiv$-broom} (Figure~\ref{fig} (a)) with respect to $y_0y_1$ and $\phiv\in \CC^\D(G-y_0y_1)$ is a sequence $B = (y_0,e_1,y_1,\ldots,e_p,y_p)$ with $p \ge 3$ such that $e_1=y_0y_1$, $e_2=y_1y_2$, $\phiv(e_2)\in \pbar(y_0)$ and for all $i \ge 3$, $e_i =y_2y_i$ and $\phiv(e_i)\in \pbar(y_j)$ for some $j<i$.

\begin{lemma}{(Cao, Chen, Jing, Stiebitz and Toft~\cite{survey})}
\label{gbroom}
Let $G$ be a $\Delta$-critical graph, $y_0y_1=e_1\in E(G)$ and $\phiv\in \mathcal{C}^\D(G-e_1)$. If $B= (y_0,e_1,y_1,\ldots,e_p,y_p)$ is a $\phiv$-broom and $\min\{d(y_1),d(y_2)\}<\D$, then the vertex set of $B$ is elementary.
\end{lemma}

\noindent
A {\it kite} $H$ (Figure~\ref{fig} (b)) is a graph with 
	$$V(H)=\{a,b,c,u,s_1,s_2,t_1,t_2\}~~  \text{and}  ~~E(H)=\{ab,ac,bu,cu,us_1,us_2, s_1t_1,s_2t_2\}.$$ 
	The lemma below reveals some properties of a kite with specified colors on its edges. 
	
		\begin{lemma}\label{kite}{(Cao, Chen and Shan~\cite{cao2020deltacritical})}
		Let $G$ be a $\Delta$-critical graph, $H\subseteq G$ 
		be a kite with $V(H)=\{a,b,c,u,s_1,s_2,t_1,t_2\}$, and let $\varphi\in \CC^\Delta(G-ab)$. 
		Suppose  that  both $K=abus_1t_1$ and $K^*=bacus_2t_2$
		are  Kierstead paths with respect to $ab$ and $\varphi$. 
		If $\varphi(s_1t_1)=\varphi(s_2t_2)$, 
		then $|\pbar(t_1)\cap \pbar(t_2) \cap ( \pbar(a)\cup \pbar(b))|\le 4$.  
	\end{lemma}

	Let $G$ be a $\Delta$-critical graph, $ab\in E(G)$, and $\varphi \in \CC^\Delta(G-ab)$. 
	A {\it fork}  $H$ (Figure~\ref{fig} (c)) with respect to $\varphi$ is a graph with 
	$V(H)=\{a,b,u,s_1,s_2,t_1,t_2\}$ \text{and}  $E(H)=\{ab,bu,us_1,us_2, s_1t_1,s_2t_2\}$
	such that $\varphi(bu)\in \pbar(a)$, $\varphi(us_1), \varphi(us_2) \in \pbar(a)\cup \pbar(b)$,  and $\varphi(s_1t_1)\in (\pbar(a)\cup \pbar(b))\cap \pbar(t_2) $
	and $\varphi(s_2t_2)\in (\pbar(a)\cup \pbar(b))\cap \pbar(t_1)$. 
	Forks  may not exist in a $\Delta$-critical graph if 
	the degree sum of $a$, $t_1$ and $t_2$ is small.

	\begin{lemma}{(Cao and Chen \cite{Cao2020-2})}\label{fork}
	Let $G$ be a $\Delta$-critical graph, $ab\in E(G)$, and $\{u,s_1,s_2, t_1,t_2\}\subseteq V(G)$. 
	If $\Delta\ge d_G(a)+d_G(t_1)+d_G(t_2)+1$, 
	then for any $\varphi\in \CC^\Delta(G-ab)$, $G$ does not contain a fork on $\{a,b,u,s_1,s_2,t_1,t_2\}$ with respect to $\varphi$.   
\end{lemma}

\subsection{Proof of Lemma~\ref{277}}

\begin{LEM2.7}
Let $G$ be a  $\Delta$-critical graph with  $\Delta\ge 7$. Then $|N_{\D\sim 2}(v)|\le 5$ for every $v\in V(G)$.
\end{LEM2.7}
\vspace{-0.3cm}
\begin{proof}
Suppose to the contrary that there is a $\D$-vertex $v$ with  $|N_{\D\sim 2}(v)|\ge 6$. By Lemma~\ref{Le:delta+2},  $v$  has no $2$-neighbors  and by Lemma~\ref{Le: VAL},  each vertex $z\in N_{\D\sim 2}(v)$ has exactly one $2$-neighbor. Let $N^2(v)=N(N(v))\backslash N[v]$. Since $|N_{\D\sim 2}(v)|\ge 6$, there are at least three 2-vertices in $N^2(v)$. 
Let $x$ be a 2-vertex in $N^2(v)$ and $y$ be a vertex in $N(x)\cap N(v)$. Clearly $y\in N_{\D\sim 2}(v)$. Let $\phiv\in \CC^\D(G-xy)$.  Then $\pbar(x)\cup\pbar(y) = C$. We first point out one  fact that will be used very often.

\noindent {\bf Fact 1}. Let  $t_1,t_2$ be two 2-vertices in $N^2(v)\backslash \{x\}$, $s_1\in N(v)\cap N(t_1)$ and $s_2\in N(v)\cap N(t_2)$. 

 (a) If $|N(x)\cap N(v)|=2$ and $\phiv(s_1 t_1)=\phiv(s_2 t_2)$, then $\phiv(t_1)\not=\phiv(t_2)$.
 
  (b) If $\phiv(s_1 t_1)\not=\phiv(s_2 t_2)$, then either $\phiv(s_1 t_1)\in \phiv(t_2)$ or $\phiv(s_2 t_2)\in \phiv(t_1)$.

\begin{proof}
 (a)  Denote  $N(x)\cap N(v)=\{y,z\}$. Suppose to the contrary that $\phiv(t_1)=\phiv(t_2)$. Then $|\pbar(t_1)\cap \pbar(t_2)|\ge 5$ since $\D\ge 7$, and $\{x,y,z,v,s_1,s_2,t_1,t_2\}$ form a kite with $\phiv(s_1 t_1)=\phiv(s_2 t_2)$, a contradiction to Lemma~\ref{kite}. 

(b) Suppose to the contrary that  $\phiv(s_1 t_1)\in \pbar(t_2)$ and $\phiv(s_2 t_2)\in \pbar(t_1)$.  Then $\{x,y,v,s_1,s_2,t_1,t_2\}$ form a fork with $\D\ge 7=d(x)+d(t_1)+d(t_2)+1$, a contradiction to Lemma~\ref{fork}.
\end{proof}

We consider two cases in  the following: there are three 2-vertices in $N^2(v)$, or there are at least four 2-vertices in $N^2(v)$.

\medskip \noindent
{\bf Case 1}: There are exactly three 2-vertices in $N^2(v)$.

Let $t_1, t_2$ be the two 2-vertices in $N^2(v)\backslash \{x\}$. Since $N_{\D\sim 2}(v)\ge 6$, we have $|N(t_i)\cap N(v)|= 2$ for each $i=1,2$ and $|N(x)\cap N(v)|=2$. Let $N(t_i)\cap N(v)=\{s_i,s_i'\}$ for each $i=1,2$. By the symmetry between $s_i$ and $s_i'$, we may assume that $\phiv(s_1 t_1)\not=\phiv(s_2 t_2)$. By Fact 1(b), we may assume $\phiv(s_1' t_1)=\phiv(s_2 t_2)$.  Applying Fact 1(a) on $s_1', t_1, s_2, t_2$,  we have $\phiv(t_1) \not = \phiv(t_2)$.  Thus $\phiv(s_2' t_2)\not=\phiv(s_1 t_1)$, $\phiv(s_2't_2) \not \in \phiv(t_1)$ and $\phiv(s_1t_1) \notin \phiv(t_2)$. This gives a contradiction to Fact 1(b) on $s_1,t_1,s_2',t_2$.

\medskip \noindent
{\bf Case 2}: There are at least four 2-vertices in $N^2(v)$.

Let $t_1,t_2,t_3$ be three 2-vertices in $N^2(v)\backslash \{x\}$, $s_i$ be a vertex in $N(t_i)\cap N(v)$, and $s_i'$ be the other neighbor of $t_i$ for each $i=1,2,3$.  

\medskip \noindent
{\bf Claim A.} $\phiv(s_i t_i) \not = \phiv(s_j t_j) $ for any $1\leq i < j \leq 3$.

\begin{proof}
Prove by contradiction.   Since $\D\ge 7> d(t_1)+d(t_2)+d(t_3)$, let $\eta\in \pbar(t_1)\cap \pbar(t_2) \cap \pbar(t_3)$.  By symmetry, we only need to consider the following two cases: $\phiv(s_1 t_1)=\phiv(s_2 t_2)=\phiv(s_3 t_3)=\alpha$, or $\phiv(s_1 t_1)=\phiv(s_2 t_2)=\alpha$ and $\phiv(s_3 t_3)=\beta\not=\alpha$.

Suppose that $\phiv(s_1 t_1)=\phiv(s_2 t_2)=\phiv(s_3 t_3)=\alpha$.  Then by symmetry, we may assume that $P_{t_1}(\alpha,\eta,\phiv)$  does not pass through  $t_2,t_3$. Let $\phiv'=\phiv/P_{t_1}(\alpha,\eta,\phiv)$. Then $s_1,t_1,s_2,t_2$ give a contradiction to Fact 1(b)  under $\phiv'$.

Suppose that $\phiv(s_1t_1)=\phiv(s_2 t_2)=\alpha$ and $\phiv(s_3 t_3)=\beta\not=\alpha$.  If $P_{t_1}(\alpha,\eta,\phiv)$ does not end at $t_2$,  let $\phiv'=\phiv/P_{t_1}(\alpha,\eta,\phiv)$. Then $s_1,t_1,s_2,t_2$ give a contradiction to Fact 1 (b) under $\phiv'$. Thus $P_{t_1}(\alpha,\eta,\phiv)$ ends at $t_2$, so $P_{t_3}(\alpha,\eta,\phiv)$ does not pass through $t_1,t_2$. Let $\phiv_1=\phiv/P_{t_3}(\alpha,\eta,\phiv)$. Now $\alpha\in \pbar_1(t_3)$. Then by Fact 1(b), we have $\phiv_1(t_1)=\phiv_1(t_2)=\{\alpha,\beta\}$. Let $\eta'\in \pbar_1(t_1)\cap \pbar_1(t_2)\cap \pbar(t_3)$. By symmetry, we may assume that $P_{t_3}(\beta,\eta',\phiv_1)$ does not pass through $t_1$. Let $\phiv_2=\phiv_1/P_{t_3}(\beta,\eta',\phiv_1)$. Then $s_1,t_1,s_3,t_3$ give a contradiction to  Fact 1(b) under $\phiv_2$.  This proves Claim A.
\end{proof}

Let $\phiv(s_1 t_1)=\alpha,\phiv(s_2 t_2)=\beta,\phiv(s_3 t_3)=\gamma$.

\medskip \noindent
{\bf Claim B.} $ \{\phiv(s_1't_1), \phiv(s_2't_2), \phiv(s_3't_3)\} =\{\phiv(s_1 t_1), \phiv(s_2 t_2), \phiv(s_3 t_3)\}$.

\begin{proof}
By  Claim A, $\alpha,\beta,\gamma$ are distinct. Suppose that $\phiv(t_1)=\{\alpha,\eta\}$ where $\eta\notin \{\beta,\gamma\}$. By Fact 1(b), we have $\phiv(t_2)=\{\beta,\alpha\}$ and $\phiv(t_3)=\{\gamma,\alpha\}$. Then $s_2,t_2,s_3,t_3$ give a contradiction to Fact 1(b). Thus by symmetry, we may assume that $\phiv(t_1)=\{\alpha,\beta\}$. Now by applying Fact 1(b) on $s_1,t_1,s_3,t_3$, we have $\phiv(t_3)=\{\alpha,\gamma\}$; By applying Fact 1(b) on $s_2,t_2,s_3,t_3$, we have $\phiv(t_2)=\{\beta,\gamma\}$. This proves Claim B.
\end{proof}

{\bf The final step.} Without loss of generality, assume $\phiv(t_1)=\{\alpha, \beta \}$. Since $|N_{\D\sim 2}|\ge 6$, let $s_4\in N_{\D\sim 2}\backslash \{s_1,s_2,s_3\}$ and  $t_4$ be the $2$-neighbor of $s_4$.  If  $t_4\in \{t_1,t_2,t_3\}$, then by symmetry, we may assume that $t_4=t_1$. Then $\phiv(s_4 t_1)=\beta$ and $s_4,t_1,s_3,t_3$ give a contradiction to Fact 1(b). If $t_4\not\in \{t_1,t_2,t_3\}$, then by Claim A, $\phiv(s_4t_4) \not = \phiv(s_it_i)$ for each $i =1,2,3$. Thus  $\{s_1,s_2,s_4,t_1,t_2,t_4\}$  does not  satisfy   Claim B. This completes the proof of  Case 2 and thus of  Lemma~\ref{277}. 
\end{proof}

\vspace{-0.3 cm}
\subsection {Proof of Lemma~\ref{Le:36}}

\vspace{-0.3 cm}

\begin{LEM2.8}
Let $G$ be a $\Delta$-critical graph and $xyrst$ be a path with $d(x)+d(y)=\Delta+2$ and $\max\{d(x),d(y)\}< \Delta$. Then  $d(t)\ge \Delta-2$.
\end{LEM2.8}

\vspace{-0.3 cm}

\begin{proof}
Let $\phiv\in \CC^\D(G-xy)$. Since  $d(x)+d(y)=\D+2$,  we have $\pbar(x)\cup\pbar(y)=C$. Let $\varphi(yr)=\alpha, \varphi(rs)=\beta, \varphi(st)=\gamma$.  Then $\alpha\in\pbar(x)$ and $\beta,\gamma\in\pbar(x)\cup\pbar(y)$.  Since $d(x) < \D$ and $d(y) < \D$, we have $|\pbar(x)|\geq 2$ and $|\pbar(y)|\geq 2$. Suppose to the contrary that $d(t)\le \D-3$. Then $|\pbar(t)|\ge 3$. 

\medskip \noindent
{\bf Claim A}.
There is a coloring  in $\mathcal{C}^\Delta(G-xy)$ such that $yr$ and $st$ are colored differently, i.e., we may assume $\alpha\not=\gamma$.

\begin{proof} 
Suppose to the contrary that $\alpha=\gamma$. Since $d(t)\le \D-3$,  let  $\eta\in  \pbar(t) \setminus \{\alpha, \beta\}$. 

If  $\eta\in \pbar(y)$, then $P_{x}(\alpha, \eta, \varphi)=P_{y}(\alpha, \eta, \varphi)$ by Lemma~\ref{vf} and thus is disjoint from  $P_{t}(\alpha, \eta, \varphi)$. Let $\varphi_1=\varphi/P_{t}(\alpha, \eta, \varphi)$. Then   $\varphi_1(yr)\not=\varphi_1(st)$, as desired.

Suppose  $\eta\in \pbar(x)$. Since $|\pbar(y)| \geq 2$, let  $\delta\in \pbar(y)\setminus \{\beta\}$.  Clearly $\delta\not \in \{\phiv(yr), \phiv(rs), \phiv(st)\}$.   Let $\varphi_1=\varphi/P_{x}(\delta, \eta, \varphi)$ and  we are back to the case when $\eta\in \pbar(y)$. This proves Claim A. 
\end{proof}

\setlength\parindent{1.7em} From now on, we assume that $\alpha\not=\gamma$ in the following proof. 

\medskip \noindent
{\bf Claim B}. We may further assume that $\alpha,\beta\in \pbar(t)$.

\begin{proof} We consider two cases: $\beta\in \pbar(t)$ and $\beta\notin \pbar(t)$.

{\bf Case B.1}: $\beta\in \pbar(t)$.

We may assume $\alpha \in \phiv(t)$ otherwise we  are done.
Let $\eta\in  \pbar(t) \setminus \{\alpha, \beta\}$. Clearly $\eta\not=\gamma$ since $\phiv(st)=\gamma$. 

If  $\eta\in \pbar(y)$,  let $\varphi_1=\varphi/P_{t}(\alpha, \eta, \varphi)$. Then  we have $\alpha,\beta\in \pbar_1(t)$, as desired. 

If   $\eta\in \pbar(x)$, let $\delta \in \pbar(y)\setminus \{\beta\}$. By Lemma~\ref{vf}, regardless of whether $\delta=\gamma$ or not, $P_{x}(\delta, \eta, \varphi)$ does not contain $yr,rs$ or $st$ since $\eta\in \pbar(t)$. Let $\varphi_1=\varphi/P_{x}(\delta, \eta, \varphi)$ and  we are back to the case when $\eta\in \pbar(y)$. This completes the proof of Case B.1.

\setlength\parindent{1.7em} {\bf Case B.2}: $\beta\notin \pbar(t)$.

\setlength\parindent{2.5em}\textbf{Case B.2.1}: $\alpha\in \pbar(t)$.

If  $\beta\in\pbar(y)$, then by Lemma~\ref{vf},  $P_x(\alpha,\beta,\phiv)$  is disjoint from $P_t(\alpha,\beta,\phiv)$. Thus $P_t(\alpha,\beta,\phiv)$  does not contain $yr$  or $rs$.  Let $\phiv_1=\phiv/P_t(\alpha,\beta,\phiv)$. Then $\beta\in \pbar_1(t)$ and we are back to  Case B.1. 

Now assume  $\beta\in \pbar(x)$.   If there is a color $\delta \in \pbar(y)\cap\pbar(t)$, let $\phiv_1 = \phiv$. Otherwise, let $\delta \in \pbar(y)$ and $\eta \in \pbar(t)\setminus\{\alpha\}$.  Then $P_t(\eta,\delta, \phiv)$  does not pass through $x$ or $y$. Let $\phiv_1 = \phiv/P_t(\eta,\delta, \phiv)$. Then $\delta \in \pbar_1(y)\cap \pbar_1(t)$ and $\beta \in \pbar_1(x)$.  Note that $P_x(\delta, \beta, \phiv_1)$ and $P_t(\delta, \beta, \phiv_1)$ are disjoint.  If $P_t(\delta, \beta, \phiv_1)$ does not contain $rs$, let $\phi_2 = \phiv_1/P_t(\delta, \beta, \phiv_1)$ and then $\phiv_2$ is a desired coloring.  If $P_x(\delta, \beta, \phiv_1)$ does not contain $rs$, let $\phi_2 = \phiv_1/P_x(\delta, \beta, \phiv_1)$. Then $\beta \in \pbar_2(y)$ and we are back to the case when $\beta\in\pbar(y)$.  This proves Case B.2.1.

\setlength\parindent{2.5em}\textbf{Case B.2.2}: $\alpha\notin \pbar(t)$.

If there is a color $\delta \in  \pbar(y)\cap \pbar(t)$,  let  $\phiv_1=\phiv/P_{t}(\alpha, \delta, \phiv)$. Then   $\alpha\in \pbar_1(t)$  and we are back to Case B.2.1.

Suppose $\pbar(y)\cap \pbar(t) = \emptyset$.  Let $\eta\in \pbar(t)$ and $\delta\in \pbar(y)\setminus \{\beta\}$. Then $\delta \in \phiv(t)$. By Lemma~\ref{vf}, regardless of whether $\delta=\gamma$ or not, $P_{x}(\delta, \eta, \phiv)$ does not contain $yr,rs$ or $st$ since $\eta\in \pbar(t)$. Let $\phiv_1=\varphi/P_{x}(\delta, \eta, \varphi)$, we are back to the case when $ \pbar(y)\cap \pbar(t) \not = \emptyset$. This completes the proof of  Case B.2 and thus the proof of  Claim B. 
\end{proof}

\setlength\parindent{1.7em} By Claim B, we assume that $\alpha,\beta\in \pbar(t)$ in the following proof. 

\medskip \noindent
{\bf Claim C}. We may further assume that $\beta,\gamma\in \pbar(y)$.

\begin{proof} We consider two cases: $\beta\in \pbar(y)$ and $\beta\notin \pbar(y)$.

{\bf Case C.1}: $\beta\in \pbar(y)$.

We may assume $\gamma\in \phiv(y)$ otherwise  we are done.  Let $\eta\in  \pbar(t) \setminus \{\alpha, \beta\}$. 

Similar to the argument in Case B.2, we may assume that there is a color $\delta \in \pbar(y)\cap\pbar(t)$ and $\delta \not = \beta$.  Then $P_t(\delta, \gamma, \phiv)$ and $P_x(\delta,\gamma, \phiv)$ are disjoint. Let $\varphi_1=\varphi/P_{x}(\delta,\gamma, \varphi)$. Then we have  $\beta,\gamma\in \pbar_1(y)$, as desired. This completes the proof of Case C.1.

\setlength\parindent{1.7em} {\bf Case C.2}: $\beta\notin \pbar(y)$.

If  $\gamma\in \pbar(y)$, then $P_t(\gamma,\beta,\phiv)$ and $P_x(\gamma,\beta,\phiv)$ are disjoint  by Lemma~\ref{vf}. Note that $rs$ and $st$ are contained in $P_t(\gamma,\beta,\phiv)$. Let $\phiv_1=\phiv/P_x(\gamma,\beta,\phiv)$. Then $\beta\in \pbar_1(y)$ and we are back to Case C.1. 

Suppose $\gamma\in \pbar(x)$.  Similar to the argument in Case B.2,  we can assume that there is a color $\delta\in \pbar(y)\cap \pbar(t)$. Then $\delta \not \in \{ \alpha,\beta\}$. Thus  $P_x(\eta,\gamma,\phiv)$ is disjoint from $P_t(\eta,\gamma,\phiv)$, so it does not contain $st$ since $\eta\in \pbar(t)$. Let $\phiv_1=\phiv/P_x(\eta,\gamma,\phiv)$ and  we are back to the case when $\gamma\in \pbar(y)$. This completes the proof of Case C.2, and thus Claim C holds.\end{proof}

Now by Claims A, B and  C,  we assume that $\phiv\in \CC^\D(G-xy)$   satisfies the following properties: 

$\bullet$ $\phiv(yr)=\alpha, \phiv(rs)=\beta, \phiv(st)=\gamma$,

$\bullet$ $\alpha \not = \gamma$,

$\bullet$ $\alpha, \beta \in \pbar(t)$ and $\beta, \gamma \in \pbar(y)$.

 Let $\phiv_1=\phiv/P_t(\alpha,\gamma,\phiv)$. Under the coloring $\phiv_1$, $P_y(\beta,\alpha,\phiv_1)=yrst$ ends at $t$ but not $x$, a contradiction to Lemma~\ref{vf}. This completes the proof of Lemma~\ref{Le:36}. \end{proof}

\vspace{-0.3 cm}
\subsection {Proof of Lemma~\ref{Le:37}}
\vspace{-0.3 cm}

\begin{LEM2.9}
Let $G$ be a $\Delta$-critical graph and $xyrst$ be a path with  $d(x)=3$ and $d(y)=\Delta$.   Suppose that $y$  has a neighbor $z \not \in \{x,r,s\}$ with  $d(z)\le \Delta-2$. Then   $d(s)\ge \Delta-1$; and $d(z)+d(t)\ge \Delta+1$ if $d(t) \le \Delta -4$.
\end{LEM2.9}
\vspace{-0.3 cm}
\begin{proof}
 Let $\phiv$ be a coloring in $\mathcal{C}^\Delta(G-xy)$.   Since $d(z)\le \Delta-2$, $d(x)=3$ and $d(y)=\Delta$,  we have $|\phiv(x)\cap \phiv(y)| = 1$. By Lemma~\ref{path}, without loss of generality, assume $\phiv(x) = \{1,2\}$, $\phiv(yz) = 2$,  $\phiv(yr) = 3$.   Denote  $\phiv(rs) = \beta$ and $\phiv(st) = \gamma$. Note that $\pbar(y)=\{1\}$.
 
  \medskip \noindent
(1) We first show $d(s) \geq \Delta -1$. 

Suppose to the contrary $d(s) \leq \Delta -2$. 

We first consider the case when  $\phiv(rs)=\beta  \not = 2= \phiv(yz)$. Then  $K=xyrs$ is a Kierstead path.  By Lemma~\ref{p4}, $|\pbar(s) \cap (\pbar(x)\cup \pbar(y)|\leq 1$. Thus $d(s) \geq 2\Delta - (d(x) + d(y) + 1 = \Delta - 2$. Since $d(s) \leq \Delta -2$, we have  $d(s) = \Delta -2$. Note that $d(s) = \Delta -2$ only if $2\in \pbar(s)$ and $|\pbar(s) \cap (\pbar(x)\cup \pbar(y))|= 1$.   Denote $\pbar(s) = \{2,\alpha\}$. 

If $ \pbar(z)\setminus \{\alpha, \beta\} \not = \emptyset$,   then $\eta  \in \pbar(z)\setminus \{\alpha,\beta\}$. By Lemma~\ref{path},  $P_z(\eta, 2, \phiv)$ ends at $x$ or $y$ and thus it does not  pass through $s$. Let $\phiv_1= \phiv/P_z(\eta, 2, \phiv)$. Then   $xyrs$  remains a Kierstead path  with respect to $\phiv_1$ and $xy$. However, $\pbar_1(s) = \{2,\alpha\} \subseteq \pbar(x)\cup  \pbar(y)$, a contradiction to Lemma~\ref{p4}.  Therefore $\pbar(z)\setminus \{\alpha, \beta\}  = \emptyset$. Since $d(z)\leq \D -2$, we have $\pbar(z) = \{\alpha, \beta\}$. 

  If $\beta \not =1$,  then  we may assume $\alpha =1$. Otherwise  both $P_z(1,\alpha, \phiv)$ and $P_s(1,\alpha, \phiv)$ are disjoint from $P_x(1,\alpha, \phiv)$. Let $\phiv_2 = \phiv/(P_z(1,\alpha, \phiv)\cup P_s(1,\alpha, \phiv))$. Then $1$ is missing at both $z$ and $s$ and $3, \beta \in \pbar_1(x)\cup \pbar_1(y)$.  Since $1 \in \pbar(z)\cap \pbar(s)$,  both $P_z(1,3, \phiv)$ and $P_s(1,3, \phiv)$ are disjoint from $P_x(1,3, \phiv)$ and thus neither  passes through $x,y$. Let $\phiv_2 = \phiv/(P_z(1,3, \phiv)\cup P_s(1,3, \phiv))$. Then $3 \in \pbar_2(z)\cap \pbar_2(s)$ and $2 \in \phiv_2(x)\cap \phiv_2(y)$.  By Lemma~\ref{path}, $P_z(2,\beta,\phiv_2)$ ends at either $x$ or $y$  and thus is disjoint from $P_s(2,\beta,\phiv_2)$. Let $\phiv_3 = \phiv_2/P_s(2,\beta,\phiv_2)$. Then $P_z(2,\beta,\phiv_3) = zyrs$ which does not end at $x$ or $y$, a contradiction to Lemma~\ref{path}.
  
 Now assume $\beta = 1$.  Then $\pbar(z) = \{1, \alpha\}$ and thus  $P_s(1,\alpha,\phiv)$ does not pass through $x,y,z$. Interchange colors on $P_s(1,\alpha, \phiv)$ and we are back to the case when $\beta \not = 1$.    Therefore this completes the proof when $\beta\not=\phiv(yz)$.
  
 Now we consider the case when  $\beta=\phiv(yz)=2$. Let $\eta$ be a color in $\pbar(z)$. Clearly $\eta\not=2$. If $\eta=1$, then by recoloring $yz$ with $1$, we are back to the case when $\beta\not=\phiv(yz)$. Thus $\eta\in \pbar(x)$. Then $P_x(\eta,1,\phiv)=P_y(\eta,1,\phiv)$. Thus by interchanging $\eta$ and $1$ on $P_x(\eta,1,\phiv)$ and then recoloring $yz$ with $\eta$, we are back to the case when $\beta\not=\phiv(yz)$. This completes the proof that $d(s) \geq \Delta -1$.
  
  \medskip \noindent
  (2) Now we assume $d(t) \leq \Delta -4$ and show  $d(z) + d(t) \geq \Delta + 1$. 
  
  Suppose to the contrary that $d(z) + d(t) \leq \Delta$.

\textbf{Claim A}.
There is a coloring in $\mathcal{C}^\Delta(G-xy)$ such that $yr$ and $st$ receive distinct colors, i.e.,  we may assume that $\gamma \not = 3$.

\begin{proof} 
Suppose to the contrary that $\gamma = 3$. Let $\eta \in  \pbar(t)\setminus \{2, 3, \beta\}$. Then $\eta \in \pbar(x)\cap \pbar(y)$.

If $\eta =1$, then  $P_{x}(3, \eta, \varphi)=P_{y}(3, \eta, \varphi)$ by Lemma~\ref{vf}, so $P_{t}(3, \eta, \varphi)$ is disjoint from $P_{x}(3, \eta, \varphi)$. Let $\varphi_1=\varphi/P_{t}(3, \eta, \varphi)$. We have that $\varphi_1(yr)\not=\varphi_1(st)$ now. 

If $\eta \not = 1$, then $P_t(1,\eta,\phiv)$ does not contain $x$ or $y$. Let $\varphi_1=\varphi/P_{t}(1, \eta, \varphi)$ and we are back to the previous case.  This proves Claim A. 
\end{proof}

\setlength\parindent{1.7em} From now on, we assume that $\phiv(yr) \not = \phiv(st)$ (i.e. $\gamma\not=3$) in the following proof. 

\medskip \noindent
{\bf Claim B}. We may further assume that $3,\beta \in \pbar(t)$.

\begin{proof} We split the proof into two cases: $\beta\in \pbar(t)$ and $\beta\notin \pbar(t)$.

{\bf Case B.1}: $\phiv(rs)=\beta\in \pbar(t)$. 

\setlength\parindent{2.5em}\textbf{Case B.1.1}:  $\beta  \not \in \pbar(y)$. Then $\beta \not = 1$.

 If $1\in \pbar(t)$, then $P_t(1,3,\phiv)$ is disjoint from $P_x(1,3,\phiv) = P_y(1,3,\phiv)$ and $yr,rs \not \in P_t(1,3, \phiv)$. Let $\phiv_1 = \phiv/P_t(1,3,\phiv)$. Then $\phiv_1(yr) = 3$, $\phiv_1(rs)= \beta$, $\phiv_1(st)=\gamma$ and $3,\beta \in \pbar_1(t)$, as desired.
 
 Now assume $1 \not \in \pbar(t)$.  Since $d(t) \leq \Delta -4$, let $\eta \in \pbar(t) \setminus \{2,3, \beta\}$. Then $\eta \in \pbar(x)$.  Thus $P_t(1,\eta,\phiv)$ does not pass through $x$ or $y$ and does not contain $yr, rs$, or $st$. Let $\phiv_1 = \phiv/P_t(1,\eta,\phiv)$ and we are back to the case when $1 \in \pbar(t)$.
 
 \setlength\parindent{2.5em}\textbf{Case B.1.2}:  $\beta  \in \pbar(y) $. Then $\beta = 1$.

If $\gamma \not = 2$, then $\gamma\in \pbar(x)$. Thus, $P_t(1,\gamma,\phiv)$ is disjoint from $P_x(1,\gamma, \phiv)$. Let $\phiv_1 = \phiv/P_t(1,\gamma,\phiv)$. Then $\phiv_1(rs) = \gamma \in \pbar_1(x)\cap \pbar_1(t)$ and $\gamma \not = 1$. We are back to Case B.1.1.

 Now assume $\phiv(st) = \gamma = 2$. Since $d(z) + d(t) \leq \Delta$ and $2 \in \phiv(z)\cap \phiv(t)$, there is a color $\eta \in \pbar(z)\cap \pbar(t)$. Since $\eta \not= \gamma$,  we have $\eta \in \pbar(x) \cup \pbar(y)$. 
 
 If $\eta \not =1$,  by Lemma~\ref{path}, $P_z(2,\eta,\phiv)$ ends at $x$.   Thus $P_t(2,\eta,\phiv)$ does not pass through $x$ or $y$ and does not contain the edge $rs$. Let $\phiv_1 = \phiv/P_t(2,\eta,\phiv)$. Then $\phiv_1(st) = \eta \in \pbar_1(x)\cup \pbar_1(y)$ and we are back to the previous case

If $\eta = 1$,  then $P_z(1,2,\phiv) = yz$. Let $\phiv_1 = \phiv/P_z(1,2,\phiv)$ and  we are back to Case B.1.1. This completes the proof of Case B.1.

\setlength\parindent{1.7em} {\bf Case B.2}: $\phiv(rs)=\beta\notin \pbar(t)$.

\setlength\parindent{2.5em}\textbf{Case B.2.1}: $\phiv(yr)=3 \in \pbar(t)$.

\setlength\parindent{3.3em}\textbf{Case B.2.1.1}: $\beta \in \pbar(y)$. That is $\beta = 1$.

Then $P_x(3,\beta,\phiv)$ ends at $y$ by Lemma~\ref{vf} and it contains both $yr$ and $rs$. Thus $P_x(3,\beta,\phiv)$ and $P_t(3,\beta,\phiv)$  are disjoint.  Let $\phiv_1=\phiv/P_t(3,\beta,\phiv)$. Then $\beta\in \pbar_1(t)$ and we are back to Case B.1.

\setlength\parindent{3.3em}\textbf{Case B.2.1.2}: $\beta = \phiv(yz)=2$.

If $1  \in \pbar(z)$, recolor $yz$ with $1$. We are back to Case B.2.1.1. 

Assume $1 \not \in \pbar(z)$. Since $d(z) \leq \Delta -2$ and $2 \in \phiv(z)$, let $\eta \in \pbar(z)\backslash\{3\}$. Clearly $\eta \not = 1,2$ and $\eta\in \pbar(x)$.  Then $P_z(1,\eta, \phiv)$ does not pass through $x$ or $y$ and does not contain the edge $rs$. Let $\phiv_1 = \phiv/P_z(1,\eta, \phiv)$. Then $1\in \pbar_1(z)$ and we are back to the previous case.

\setlength\parindent{3.3em} \textbf{Case B.2.1.3}: $\beta  \in \pbar(x)$.  

We may further assume $1\in \pbar(t)$. Otherwise, since $d(t)\le \D-4$, let $\eta\in  \pbar(t) \setminus \{ 2,3\}$.  Then $\eta \not \in \{1,2,3, \beta\}$,  and  $P_x(1,\eta,\phiv)$ and $P_t(1,\eta, \phiv)$ are disjoint.  Let $\phiv_1 = \phiv/P_t(1,\eta,\phiv)$. Then $1\in \pbar_1(t)$.

  Note $P_x(\beta,1,\phiv_1)$ and $P_t(\beta,1,\phiv_1)$  are disjoint. If $P_x(\beta,1,\phiv_1)$ does not contain the edge $rs$, let $\phiv_2 = \phiv_1/P_x(\beta,1,\phiv_1)$ and we are back  to Case B.2.1.1. If $P_t(\beta,1,\phiv_1)$ does not contain the edge $rs$, let $\phiv_2 = \phiv_1/P_t(\beta,1,\phiv_1)$ and we are back  to Case B.1.  This completes the proof of Case B.2.1.

\setlength\parindent{2.5em}\textbf{Case B.2.2}: $3\notin \pbar(t)$.

Since $d(t)\le \D-4$, let  $\eta\in  \pbar(t) \setminus \{2, 3, \beta\}$. 

If  $\eta=1 $, then $P_{x}(3, 1, \varphi)$ and  $P_{t}(3,1, \varphi)$ are disjoint. Let $\varphi_1=\varphi/P_{t}(1,3, \varphi)$. Then  $\phiv_1(yr) = 3 \in \pbar_1(t)$  and we are back to Case B.2.1.

Therefore $\eta\not=1$. If $\beta\not=1$, then $P_{x}(1, \eta, \varphi)$ does not contain $yr,rs$ or $st$ since $\eta\in \pbar(t)$. Let $\varphi_1=\varphi/P_{x}(1, \eta, \varphi)$ and  we are back to the case when $\eta = 1$.

If  $\beta=1$, then $P_{x}(\eta, 1, \varphi)$ and $P_{t}(\eta, 1, \varphi)$ are disjoint. If  $P_{x}(\eta, 1, \varphi)$ does not pass through $rs$, let $\varphi_1=\varphi/P_{x}(\eta, 1, \varphi)$. Then $\eta$ is missing at $y_1$ now and we are back to the case when $\eta = 1$.  If $P_{t}(\eta, 1, \varphi)$ does not contain $rs$, let $\varphi_1=\varphi/P_{t}(\eta, 1, \varphi)$. Then $\beta\in \pbar_1(t)$ and we are back to Case B.1. This completes the proof of Case B.2, and so Claim B holds.
\end{proof}

\setlength\parindent{1.7em} By  Claims A and B, we assume  that  $\phiv$ satisfies the following properties: 

$\bullet$  $\phiv(yr) =3\in \pbar(t), \phiv(rs) =\beta\in \pbar(t)$.

$\bullet$ $\phiv(st) = \gamma \not= 3$  

\medskip \noindent
{\bf Claim C}. We may further assume  $\beta = \phiv(yz) = 2$.

\begin{proof}
Suppose to the contrary  $\beta \not = 2$.

{\bf Case C.1}: $\gamma \not =\phiv(yz)$ (i.e. $\gamma\not=2$).

\setlength\parindent{2.5em}\textbf{Case C.1.1}: $1\in \{\gamma,\beta\}$.

If $\beta = 1$, then $P_t(\gamma, 1, \phiv)$  does not  pass through $x$ or $y$. Let $\phiv_1 = \phiv/ P_t(\gamma, 1, \phiv)$. Then $\phiv_1(st) = 1$.    Thus we assume  $\gamma =1$.

If  $\beta\in \pbar(z)$,  let $\phiv_1 = \phiv/P_x(\beta,1,\phiv)$ and then recolor $yz$ with $\beta$.   Then $\phiv_1$ is a desired coloring.

If  $3\in \pbar(z)$, let $\phiv_1=\phiv/P_x(\beta,1,\phiv)$ and  $\phiv_2=\phiv_1/P_z(3,\beta,\phiv_1)$. Notice that the second Kempe exchange will not effect $yr$ or $rs$ since they are on $P_x(3,\beta,\phiv_1)=P_y(3,\beta,\phiv_1)$ by Lemma~\ref{vf}. Thus we obtain a desired coloring  by recoloring $yz$ with $\beta$ under $\phiv_2$.  

Now we  assume $3, \beta \not \in \pbar(z)$.

 If $\pbar(z)\cap \pbar(t) \not =\emptyset$, let  $\eta \in \pbar(z)\cap \pbar(t)$. Then $\eta \not \in \{1,2,3, \beta\}$ and  $\eta\in \pbar(x)$. Note that $P_x(1,\eta,\phiv)=P_y(1,\eta,\phiv)$ does not contain $st$ since $\eta\in\pbar(t)$.  Let $\phiv_1=\phiv/P_x(1,\eta,\phiv)$ and then  $\eta\in \pbar_1(y)$.  Let $\phiv_2=\phiv_1/P_z(\eta,3,\phiv_1)$ and then  $3\in \pbar_2(z)$. Note that $P_z(\eta, 3, \phiv_1)$ does not contain $yr$ or $t$ since $yr$ is on $P_x(\eta,3,\phiv_1)=P_y(\eta,3,\phiv_1)$ and $3,\eta\in\pbar_1(t)$. Finally let $\phiv_3=\phiv_2/P_x(\eta,1,\phiv_2)$. We are back to the case when $\phiv(yr) \in \pbar(z)$. 

Now assume  $\pbar(z)\cap \pbar(t) =\emptyset$.   Since $d(z) + d(t) \leq \Delta$, $\phiv(z)$ and $\phiv(t)$ form a partition of $C$.  Consequently, we have $1\in \pbar(z)$ and $2 \in \pbar(t)$.   Since  $d(z) \leq \Delta -2$, let  $\eta\in \pbar(z)\setminus\{1\}$.   Clearly $\eta\in \pbar(x)$ and $\eta\notin \{1,2,3,\beta\}$. Let $\phiv_1$ be the coloring obtained from $\phiv$ by recoloring $yz$ with $1$. Then $2\in \pbar_1(y)\cap\pbar_1(z)$ and $P_x(2,\eta,\phiv_1)=P_y(2,\eta,\phiv_1)$ by Lemma~\ref{vf}.  Let $\phiv_2=\phiv_1/P_x(2, \eta,\phiv_1)$ and $\phiv_3$ be the coloring obtained from $\phiv_2$ by recoloring $yz$ with $\eta$. Now we have $\gamma =1 \in \pbar_3(y)$, $\phiv_3(yz)=\eta\not=\beta$ and $2 \in \pbar_3(z)\cap \pbar_3(t)$. Thus we are back to the case when $\pbar(z)\cap \pbar(t) \not =\emptyset$. This completes the proof of Case C.1.1.

\textbf{Case C.1.2}: $1 \notin \{\gamma,\beta\}$.

 Since $d(t)\le \Delta-4$, let $\eta\in \pbar(t)\backslash \{2, 3,\beta\}$.   We may assume $\eta = 1$. Otherwise, $\eta\in \pbar(x)$ since $\pbar(x)=C\backslash \{1, 2\}$. Thus by interchanging colors on $P_t(1,\eta,\phiv)$, $1$ is missing at $t$. Since  $\gamma\in \pbar(x)$,  we have $P_x(\gamma,1,\phiv)=P_y(\gamma,1,\phiv)$. Since $1\in \pbar(t)$, $P_x(\gamma,1,\phiv)$ does not contain $st$. Therefore, by interchanging $\gamma$ and $1$ on $P_x(\gamma,1,\phiv)$, we are back to Case C.1.1. This completes the proof of Case C.1.

\setlength\parindent{1.7em}{\bf Case C.2}: $\gamma = \phiv(yz)= 2$.

In this case, $\phiv(yz)=\phiv(st)=2\in \phiv(z)\cap \phiv(t)$.  If $1\in \pbar(z)$, recolor $yz$ with $1$. Then we are back to Case C.1 if $\beta \not = 1$. Otherwise, we have a desired coloring. Thus in the following we assume $1 \in \phiv(z)$.

\setlength\parindent{2.5em}\textbf{Case C.2.1}: $ \{3,\beta\}\cap \pbar(z) \not = \emptyset$.

 If $\beta \in \pbar(z)$, then  by Lemma~\ref{path},  $P_z(2,\beta,\phiv)$  ends at $x$ since $\beta \in \pbar(x)$ and  it is disjoint from $P_t(2, \beta, \phiv)$.  Thus $\phiv_1 = \phiv/ P_z(2,\beta,\phiv)$ is a desired coloring.
 
 Assume $3 \in \pbar(z)$ and $\beta \in \phiv(z)$. 
 
 If $\beta = 1$, then $P_y(1,3,\phiv)$ contains the edges $yr$ and $rs$ and is disjoint from $P_z(1,3,\phiv)$. Note that $1, \beta\in \pbar(t)$. Let $\phiv_1= \phiv/P_z(1,3,\phiv)$ and we are back to the case when $1 \in \pbar(z)$. 
 
  Assume  $\beta \not = 1$.  Since $d(z) \leq \Delta-2$, let $\eta \in \pbar(z)\setminus \{3\}$. Then $\eta \not \in \{1, 2, 3, \beta\}$. Thus $P_z(1,\eta,\phiv)$ does not contain the vertices $x,y$ or the edges $rs, st$. Let $\phiv_1 = \phiv/P_z(1,\eta,\phiv)$ and we are back to the case when $1\in \pbar(z)$. This completes the proof of  Case C.2.1.

\textbf{Case C.2.2}: $\{3,\beta\} \cap \pbar(z) = \emptyset$.

Since $2\in \phiv(z)\cap \phiv(t)$ and $d(z)+d(t) \leq \Delta$,  let $\eta \in \pbar(t)\cap \pbar(z)$. Then $\eta\in \pbar(x)$.  If $\beta \not = 1$,  by interchanging colors on  $P_x(\eta,1,\phiv)$  and then recoloring $yz$ with $\eta$, we are back to Case C.1. Suppose $\beta = 1$. Then $P_x(\eta,1,\phiv)$ and $P_z(\eta,1,\phiv)$ are disjoint and   either $P_x(\eta,1,\phiv)$ or  $P_z(\eta,1,\phiv)$ does not contain $rs$.  In the former case, by interchanging $\eta$ and $1$ on $P_x(\eta,1,\phiv)$ and then recoloring $yz$ with $\eta$, we are back to Case C.1.  In the later case by interchanging $\eta$ and $1$ on $P_z(\eta,1,\phiv)$ and then recoloring $yz$ with $1$, we have a desired coloring. This completes the proof of Case C.2.2, and so Claim C holds.
\end{proof}

\setlength\parindent{1.7em} By Claim C, we further assume   $\phiv(yz)= \phiv(rs) = 2$. Note that $\phiv(x)\cap \phiv(y)=\{2\}$ and $\pbar(x)\cup \pbar(y)= C\backslash \{2\}$.

\medskip \noindent
{\bf Claim D}. We may further assume that $ \pbar(y)\cap \pbar(z) \not = \emptyset$ and $\gamma \in \pbar(y)\cap \pbar(z)$. That is $\gamma = 1\in\pbar(z)$.

\begin{proof}
We split the proof into the following cases.

\setlength\parindent{2.5em} {\bf Case D.1}: $\phiv(yr) = 3\in \pbar(z)$.

\setlength\parindent{3.3em} {\bf Case D.1.1}: $\gamma = 1$.

In this case $P_x(1, 3, \phiv)$ is disjoint from  $P_z(1, 3, \phiv)$.  Let $\phiv_1 = \phiv/P_z(1, 3, \phiv)$.  If $P_z(1, 3, \phiv)$ does not end at $t$, then $\phiv_1$ is a desired coloring.   If $P_z(1, 3, \phiv)$  ends at $t$,   let $\phiv_2$ be the coloring obtained from $\phiv_1$ by recoloring  $yz$ with $1$.  In the coloring $\phiv_2$,  $2$ is missing at $y$, $3$ is missing at $x$, and $P_y(3,2,\phiv_2) = yrst$,  a contradiction to Lemma~\ref{vf}. This proves Case D.1.1.

{\bf Case D.1.2}: $\gamma \not = 1$. Then $\gamma \not \in \{1,2,3\}$ and $\gamma \in \pbar(x)$.

If $1 \in \pbar(t)$, then $P_x(1,\gamma,\phiv)$ ends at $y$ and thus  does not contain  the edge $st$. Thus by interchanging $1$ and $\gamma$ on $P_x(1,\gamma,\phiv)$, we are back to Case D.1.1.

Assume $1 \not \in \pbar(t)$. Since $d(t)\le \D-4$, let $\eta\in \pbar(t)\backslash \{2, 3\}$. Then  $\eta \not \in \{1,2,3,\gamma\}$ and $\eta\in  \pbar(x)$. By interchanging the colors  on $P_t(\eta,1,\phiv)$, we are back to the case when $1\in \pbar(t)$. This proves Case D.1.

\setlength\parindent{2.5em} {\bf Case D.2}: $\phiv(yr) = 3\notin \pbar(z)$.

Since $d(z)+d(t)\le \D$, either $\phiv(z)$ and $\phiv(t)$ form a partition of $C$ or there exists a color $\eta\in \pbar(z)\cap \pbar(t)$.

\setlength\parindent{3.3em} {\bf Case D.2.1}: There exists a color $\eta\in \pbar(z)\cap \pbar(t)$.

In this case we have $\eta\notin \{2, 3,\gamma\}$  and  $\eta\in \pbar(x)\cup \pbar(y)$.

If $\eta =1$, then $P_z(1,3,\phiv)$ does not pass through $x,y$ or $t$ since both $\alpha$ and $\eta$ are missing at $t$.   We are back to Case D.1 by interchanging $1$ and $3$ on $P_z(1,3,\phiv)$.

If $\eta \not = 1$,  then $\eta\in \pbar(x)$ and $P_x(\eta,1,\phiv)$  does not pass through $t$ since $\eta\in \pbar(t)\cap \pbar(z)$. Thus by interchanging $\eta$ and $1$ on $P_x(\eta, 1,\phiv)$, we are back to the case when $\eta =1$. This completes the proof of Case D.2.1. 

{\bf Case D.2.2}: $\phiv(z)$ and $\phiv(t)$ form a partition of $C$.

In this case $\gamma\in \pbar(z)$. If $\gamma = 1$, then $\phiv$ is a desired coloring.   Therefore we assume $\gamma \not = 1$. Thus $\gamma\in \pbar(x)$. Let   $\eta\in \pbar(t)\backslash \{2, 3\}$. By Lemma~\ref{vf}, $P_x(1, \eta,\phiv)$ does not pass through $z$ or $t$.  Note that if $1= \eta$, then $P_x(1, \eta,\phiv) = x$.  Let $\phiv_1=\phiv/P_x(1, \eta,\phiv)$.  Then $P_x(\eta,\gamma,\phiv_1)=P_y(\eta,\gamma,\phiv_1)$. Note that $P_x(\eta,\gamma,\phiv_1)$ does not contain $t$ since $\eta\in \pbar_1(t)$. Let $\phiv_2=\phiv_1/P_x(\eta,\gamma,\phiv_1)$. Then we have $\gamma\in \pbar_2(y)\cap \pbar_2(z)$ and thus $\phiv_1$ is a desired coloring.  This completes the proof of Case D.2, and so Claim D holds.
\end{proof}

In summary, by Claims A, B, C, and D, we assume that $\phiv$ satisfies the following properties: 

\noindent
$\bullet$ $\phiv(x) = \{1,2\}$ and $1\in \pbar(y)\cap \pbar(z)$

\noindent
$\bullet$ $\phiv(yr) = 3$, $\phiv(yz) = \phiv(rs) = 2$, and $\phiv(st) = 1$

\noindent
$\bullet$ $2, 3 \in \pbar(t)$.

Note that $P_x(1,3,\phiv)$ ends at $y$ and is disjoint from $P_t(1,3,\phiv)$.  If $P_t(1,3,\phiv)$ does not end at $z$, let $\phiv_1$ be the coloring obtained from $\phiv$ by interchanging colors on  $P_t(1,3,\phiv)$ and recoloring $yz$ with $1$. Then $3 \in \pbar_1(x)$, $2\in \pbar_1(y)$ and $P_y(3,2,\phiv_1)= yrst$ not ending at $x$, a contradiction to Lemma~\ref{vf}. Thus $P_t(1,3,\phiv)$ ends at $z$.  Let $\phiv_2 = \phiv/P_t(1,3,\phiv)$. Then $P_z(2,3, \phiv_2) = zyrst$ which does not end at $x$, a contradiction to Lemma~\ref{path}.  This completes the proof of Lemma~\ref{Le:37}.              
\end{proof}

\vspace{-0.3 cm}
\subsection {Proof of Lemma~\ref{Le:D+3}}

\vspace{-0.3 cm}
\begin{LEM2.11}
Let $G$ be a $\Delta$-critical graph and  $xy$ be  an edge with $d(x)+d(y)=\Delta+3$ and $\max\{d(x),d(y)\}<\Delta$. Then $x$ has $d(x)-2$ neighbors of degree $\Delta$  having no $(\Delta-2)^-$-neighbors other than $x, y$.
\end{LEM2.11}
\vspace{-0.3 cm}
\begin{proof}
 Let $\phiv \in \mathcal{C}^\Delta(G-xy)$. Since $G$ is $\D$-critical and $d(x)+d(y)=\D+3$,  we have $|\phiv(x)\cap \phiv(y)|=1$. Let $\delta$ be the color in $\phiv(x)\cap \phiv(y)$. Then $\pbar(x)\cup \pbar(y)=C\backslash \{\delta\}$. By Lemma~\ref{Le: VAL}, $x$ has at least $d(x)-2$ neighbors of degree $\D$. Thus including $y$, $x$ has at most two neighbors of degree less than $\D$. By Lemmas~\ref{p4} and \ref{path}, we have the following fact which will be applied frequently.

\medskip \noindent
{\bf Fact 1}. Let $yxzt$ be a path  with  $\phiv(xz) \in \pbar(y)$.

 (1)   $\pbar(z) \subseteq \{\delta\}$ and thus $d(z) \geq \Delta -1$. If $\delta \in \pbar(z)$, then for any color $\eta \in \phiv(z)\setminus \{\phiv(xz)\}$, $P_z(\delta, \eta, \phiv)$ ends at $x$ or $y$.
 
 (2) If $yxzt$ is a Kierstead path, then  $\pbar(t) \subseteq \{\delta\}$ and thus $d(t) \geq \D-1$.

We consider two cases in the following according to the number of $\Delta$-neighbors of $x$.

\medskip \noindent
{\bf Case 1.} $x$ has a neighbor $z_0\not=y$ with $d(z_0)<\D$.

It is sufficient to show that for any path $yxzt$ with $z\not = z_0$,  we have $d(t) \geq \Delta -1$.

Suppose to the contrary that  there is a path $yxzt$ such that $z\not = z_0$ but $d(t) \leq \Delta -2$.   We consider two cases according to $\phiv(xz_0) = \delta$ or not.

{\bf Case 1.1}: $\alpha =\phiv(xz_0)\not = \delta$.

By Fact 1(1),  $\pbar(z_0) = \{\delta\}$.

First assume $\phiv(xz)\in \pbar(y)$. Then by Fact 1(2), $\phiv(zt) = \delta$ otherwise $yxzt$ is  a Kierstead path.
    Since $d(t) \leq \D-2$ and $\delta \in \phiv(t)$, let $\eta \in \pbar(t) \setminus \{\alpha\}$.  By Fact 1(1), $P_{z_0}(\delta, \eta,\phiv)$ ends at $x$ or $y$ and thus is disjoint from $P_t(\delta, \eta, \phiv)$. Let $\phiv_1 = \phiv/P_t(\delta, \eta, \phiv)$. Then  $yxzt$ is a Kierstead path  in $\phiv_1$ and thus $d(t) \geq \Delta -1$ by Fact 1(2), a contradiction.
    
Now assume $\phiv(xz)=\delta$.  Denote $\beta = \phiv(zt)$. Then $\beta \in \pbar(x)\cup \pbar(y)$.  We may assume that $\beta\in \pbar(x)$. Otherwise if there is a color $\eta \in \pbar(t)\cap \pbar(x)$, interchange colors on the path $P_t(\eta,\beta, \phiv)$ which does not contain $x$ or $y$. If no such $\eta$ exists, let $\eta \in \pbar(x)$ and $\gamma \in \pbar(t)\setminus \{\delta\}$. Let $\phiv_1= \phiv/P_t(\eta,\gamma,\phiv)$ and then let $\phiv_2 = \phiv_1/P_t(\eta,\beta, \phiv_1)$.

By Fact 1(1), $P_{z_0}(\delta, \beta, \phiv)$ ends at $x$ and thus contains $xzt$. This implies $\delta \in \phiv(t)$.  Thus $|\pbar(t)\cap (\pbar(x)\cup \pbar(y))|\geq 2$ since $d(t) \leq \Delta -2$. Let $\eta \in \pbar(t)\setminus \{\alpha\}$. By Fact 1(1) again, $P_{z_0}(\delta, \eta,\phiv)$ ends at $x$ or $y$ and thus is disjoint from $P_t(\delta,\eta,\phiv)$. Let $\phiv_1 = \phiv/P_t(\delta,\eta,\phiv)$. Then in $\phiv_1$, $P_x(\delta, \beta, \phiv_1) = xzt$ which is disjoint from $P_{z_0}(\delta, \beta, \phiv_1)$, a contradiction to Fact 1(1).  This completes the proof of Case 1.1.

{\bf Case 1.2}: $\phiv(xz_0)=\delta$.

Then $\phiv(xz)\in \pbar(y)$.   Since $d(t) \leq \Delta -2$, by Fact 1(2), $\phiv(zt) = \delta$. Let $\eta\in \pbar(z_0)$.  Similar to  the argument in Case 1.1,  we  assume $\eta\in \pbar(x)$.   Recolor $xz_0$ with $\eta$. Then $yxz_t$ is a Kierstead path. By Fact 1(2),  $d(t) \geq \Delta -1$, a contradiction. This completes the proof of Case 1.

\medskip \noindent
{\bf Case 2.} All vertices in $N(x)\backslash \{y\}$ are $\D$-vertices. 

Since $|\pbar(y)\cap \phiv(x)|= d(x)-2$,   we are done  if  $d(t) \geq \Delta -1$ for every path $yxzt$ with $\phiv(xz) \in \pbar(y)$. Thus  assume  that there is a path $yxz_0t_0$ such that $\phiv(xz_0) \in \pbar(y)$ and $d(t) \leq \Delta -2$. By Fact 1(2), $\phiv(z_0t_0) = \delta$. Denote $\alpha = \phiv(xz_0)$. Then $\alpha\in \pbar(y)$. With  a similar argument as before, we may assume $\alpha \in \pbar(t_0)$ and there is a color $\eta \in \pbar(t_0)\cap \pbar(x)$. Then $\eta \not = \alpha$. Now it is sufficient to show that for any path $yxzt$ with $z\not = z_0$,  we have $d(t) \geq \Delta -1$.  We consider the following two cases.

{\bf Case 2.1}.  $\phiv(xz)  = \beta \in \pbar(y)$. 

Then by Fact 1(2), $\phiv(zt) = \delta$,  so $t\not=t_0$.  Since $d(t) \leq \Delta -2$, there is a color $\eta_1\in \pbar(t)$. Then $\eta_1 \in \pbar(x)\cup \pbar(y)$. Similarly we may assume $\eta,\eta_1 \in \pbar(x)$. Note that $d(t)\le \D-2$ and $d(t_0) \leq \Delta -2$. Thus $\eta \not = \eta_1$ since otherwise both $P_{t_0}(\delta, \eta,\phiv)$ and $P_{t}(\delta, \eta,\phiv)$ end at $x$ by Fact 1(2), a contradiction.

Now let $\phiv_1$ be the coloring obtained from $\phiv$ by coloring $xy$ with $\alpha$, leaving $xz_0$ uncolored and recoloring $z_0t_0$ with $\alpha$. Then $P_x(\eta_1,\delta,\phiv_1)=P_{z_0}(\eta_1,\delta,\phiv_1)$ by Lemma~\ref{vf}. Let $\phiv_2=\phiv_1/P_{t}(\eta_1,\delta,\phiv_1)$. Then $\phiv_2(zt)=\eta_1\in \pbar_2(x)$. Note that the last Kempe exchange may affect the colors of the edges incident to $t_0$, so $\delta$ may not be missing at $t_0$ under $\phiv_2$. But we still have $\eta\in \pbar_2(x)\cap \pbar_2(t_0)$.  If $\delta \in \phiv_2(t_0)$, let $\phiv_3=\phiv_2/P_{t_0}(\eta,\delta,\phiv_2)$. Otherwise let $\phiv_3 = \phiv_2$.  Then we have $\delta\in \pbar_3(z_0)\cap \pbar_3(t_0)$. Finally let $\phiv_4$ be the coloring obtained from $\phiv_3$ by recoloring $z_0t_0$ with $\delta$, coloring $xz_0$ with $\alpha$ and leaving $xy$ uncolored. Then $yxzt$ is a Kierstead path under $\phiv_4$. However $d(t) \leq \Delta -2$, a contradiction to Fact 1(2).

{\bf Case 2.2} $\phiv(xz) = \delta$.

Denote $\phiv(zt) = \beta$.  With similar arguments as before we may assume that there is a color $\eta'\in \pbar(t)\cap \pbar(x)$. We may then assume that $\beta \in \pbar(x)$ since otherwise we can interchange $\beta$ and $\eta'$ on $P_x(\beta,\eta',\phiv)$ to get a desired coloring. 

Since $d(t) \leq \Delta -2$, let $\eta_1\in \pbar(t)\backslash \{\alpha\}$. We then show that we may assume $\eta_1 \in \pbar(x)\cup \{\delta\}$. Suppose otherwise $\eta_1\in \pbar(y)\backslash \{\alpha\}$. Since $d(x) \leq \D-1$,  we have $|\pbar(x)| \geq 2$. Let $\alpha'$ be a color in $\pbar(x)\backslash \{\phiv(zt)\}$. By interchanging $\eta_1$ and $\alpha'$ on $P_x(\eta_1,\alpha',\phiv)$, we obtain a coloring as desired.

Let $\phiv_1$ be the coloring obtained from $\phiv$ by coloring $xy$ with $\alpha$, leaving $xz_0$ uncolored and recoloring $z_0t_0$ with $\alpha$. Then under $\phiv_1$, $z_0xzt$ is a Kierstead path with $\eta_1\in (\pbar_1(x) \cup \pbar_1(z_0))\cap\pbar_1(t)$,  a contradiction to Lemma~\ref{p4}. This completes the  proof of the lemma.
\end{proof}

\vspace{-0.3 cm}

\subsection {Proof of Lemma~\ref{Le:5-vertex}}
\vspace{-0.3 cm}
\begin{LEM2.12}
Let $G$ be a $7$-critical graph and $x$ be a  5-vertex.

\noindent
(1) if $x$ has three 
  6-neighbors, then each $7$-neighbor of $x$ has exactly one $5^-$-neighbor. 
  
  \noindent
  (2)  if $x$  has two 6-neighbors, then   $x$ has  two $7$-neighbors, each of which  has at most two $5^-$-neighbors.
  
  \noindent
  (3) if $x$ has exactly four 7-neighbors,  then $x$ has  two $7$-neighbors, each of which  has at most three $5^-$-neighbors.
\end{LEM2.12}

\vspace{-0.2cm}
\begin{proof}  If $x$ has a $5$-neighbor, then by Lemma~\ref{Le: VAL}, $x$ has at least three $7$-neighbors and thus has at most one $6$-neighbor. To show the lemma in this case, we only need to consider the case when $x$ has four $7$-neighbors and one $5$-neighbor which is (3), and it follows from  Lemma~\ref{Le:D+3}.  In the rest of the proof, we assume that $x$ has no $5$-neighbors.  By the assumption of the lemma, $x$ has a $6$-neighbor. Let $y$ be a $6$-neighbor of $x$, $\phiv\in \mathcal{C}^\Delta(G-xy)$. Without loss of generality  we assume that $\pbar(y)=\{1,2\}$, $\pbar(x)=\{3,4,5\}$, and  $\phiv(x)\cap \phiv(y)=\{6,7\}$. By Lemma~\ref{Le: VAL}, $x$ has at least two $7$-neighbors. 

\medskip \noindent
 {\bf (1)} Denote the two 6-vertices in $N(x)\setminus \{y\}$ by $z_1,z_2$, the two 7-vertices in $N(x)$ by $v_1,v_2$. We need to show that for any path $yxvt$ with $v\in \{v_1, v_2\}$, $d(t) \leq 5$. We consider three cases.

{\bf Case 1.1 } $x,y,z_1,z_2$ form the vertex set of a multi-fan with respect to $xy$ and $\phiv$.

In this case, by Lemma~\ref{vf}, we have $\pbar(z_1)\cup\pbar(z_2)=\{6,7\}$. Assume without loss of generality that $\pbar(z_1)=\{6\}$ and $\pbar(z_2)=\{7\}$.    Then  for each $\alpha \in \pbar(x)\cup \pbar(y)$,    both $P_{z_1}(6, \alpha,\phiv)$ and $P_{z_2}(7, \alpha,\phiv)$ end at $x$ if  $\alpha \in \pbar(x)$.

Let $yxvt$ be a path where $d(v) = 7$. Let $\eta$ be a color in $\pbar(t)$ and $\beta = \phiv(vt)$. We may assume that $\eta\in \pbar(x)$ since otherwise $\eta\in \{1,2,6,7\}$, and we can interchange $\eta$ and $3$ on $P_t(\eta,3,\phiv)$, which doesn't pass through $x$ or $y$ by Lemma~\ref{vf},  to obtain a  desired coloring.  Thus we assume $\eta\in \pbar(x)$. 

 We may further assume that $\beta=\phiv(v_1t)\in \pbar(x)$. Otherwise  $\beta \in \{1,2,6,7\}$. Note that $P_{t}(\beta, \eta,\phiv)$ does not end at $x$ or $y$.  Let $\alpha \in  \pbar(x)\setminus \{\eta\}$. Interchange $\eta$ and $\phiv(vt)=\beta$ on $P_t(\beta,\eta,\phiv)$  first and then interchange $\beta,\alpha$ on the $(\beta,\alpha)$-chain starting at $t$. We obtain a desired coloring. Thus we assume that $\beta\in \pbar(x)$.
 
Now let $\phiv_1 = \phiv/P_t(\eta,\phiv(xv),\phiv)$. Then $\phiv(xv)\in \pbar_1(t)$ and $P_x(\phiv(xv),\phiv(vt),\phiv_1)=xvt$ does not end at $y$, $z_1$, or $z_2$, a contradiction to Lemma~\ref{vf}. This completes the proof of Case 1.1.

{\bf Case 1.2} $x,y,z_1,z_2$ do not form the vertex set of a multi-fan with respect to $xy$ and $\phiv$, and $|\{\phiv(xz_1),\phiv(xz_2)\}\cap\{1,2\}|=1$.

By symmetry,  assume that  $\phiv(xz_1)=1$, $\pbar(z_1)=\{6\}$,   $\phiv(xz_2)=7$, $\phiv(xv_1)=2$, and $\phiv(xv_2)=6$. Then  for each color $\eta \in \{2,3,4,5\}$, $P_{z_1}(\eta,6,\phiv)$ ends at $x$ or $y$ depending on whether $\eta \in \pbar(x)$ or $\eta \in  \pbar(y)$ by Lemma~\ref{vf}. Similar to the argument in Case 1.1, we may further assume $3 \in \pbar(z_2)$. 

Let $yxvt$ be a path where $d(v) = 7$.  Then $\phiv(xv) \in \{2,6\}$.  We first assume  $\phiv(xv) = 2$.   

If $\phiv(vt)\in \pbar(x)\cup \pbar(y)$, then $yxvt$ is  a Kierstead path with $d(x)<\D$. Thus $\pbar(t_1)=\{6,7\}$ by Lemma~\ref{p4}. Let $\eta$ be a color in $\pbar(x)\setminus  \{\phiv(vt)\}$.  Then by Lemma~\ref{p4link}, $P_{t}(\eta,6,\phiv)$ ends at $x$. However, by Lemma~\ref{vf}, $P_{z_1}(\eta,6,\phiv)$ ends at $x$, a contradiction.

If $\phiv(vt) = 7$, recolor  $xz_2$ with $3$ and we are back to the case when $\phiv(vt)\in \pbar(x)\cup \pbar(y)$. 

If $\phiv(vt) = 6$, let $\eta \in \pbar(t) \cap (\pbar(x)\cup \pbar(y))$.  We may assume $\eta \in \pbar(x)$ since otherwise  we can pick a color $\beta \in \pbar(x)$ and interchange colors on $P_{t}(\eta,\beta, \phiv)$.  Since $P_{z_1}(\eta,6,\phiv)$ ends at $x$, $P_{t}(\eta,6,\phiv)$ and $P_{z_1}(\eta,6,\phiv)$  are disjoint.  Interchange colors on $P_{t}(\eta, 6,\phiv)$ and we are back to the case when $\phiv(vt)\in \pbar(x)\cup \pbar(y)$ again.

Now we assume  $\phiv(xv)=6$. Denote $\phiv(vt) = \beta$. 
 If $\beta = 7$, then  recolor the edge $xz_2$ with $3$ and then $7$ is missing at $x$. Thus we may assume $\beta \in  \pbar(x)\cup \pbar(y)$.  

 If  $\phiv(vt) = \beta \in \pbar(x)$, then  $P_x(6,\beta,\phiv)$ ends at $z_1$ and thus $6\in \phiv(t)$. Since $d(t) \leq 5$,  let $\alpha \in  \pbar(t)\cap (\pbar(x)\cup \pbar(y))$.  Similarly as before we may  further assume that $\alpha \in \pbar(x)$.  Note that  $P_{z_1}(\alpha,6,\phiv)$ and    $P_{t}(\alpha,6,\phiv)$ are disjoint. Let $\phiv_1= \phiv/ P_{t}(\alpha,6,\phiv)$. Then $6$ is missing at $t$ and thus $P_x(6,\beta,\phiv_1)=xvt$ does not end at $z_1$,  a contradiction.
 
Suppose $\phiv(vt)=\beta\in \pbar(y)$. Let $\alpha'$ be a color in $\pbar(t)\backslash \{7\}$. Then similarly, we can assume that $\alpha'\in \pbar(x)$. By interchanging $\alpha'$ and $\beta$ on $P_{t}(\alpha',\beta,\phiv)$, we are back to the case when $\phiv(vt)\in \pbar(x)$. This completes the proof of Case 1.2.
 
{\bf Case 1.3} $\{\phiv(xz_1),\phiv(xz_2)\}=\{6,7\}$.

Let $yxvt$ be a path where $d(v) = 7$. Without loss of generality, assume  $\phiv(xz_1)=6$, $\phiv(xz_2)=7$, and $\phiv(xv) = 1$.  Denote $\phiv(vt) = \beta$.
 
We first assume  $\beta\in \pbar(x)\cup \pbar(y)$. Then $yxvt$ is  a Kierstead path with $d(x)<\D$. Thus $\pbar(t)=\{6,7\}$ by Lemma~\ref{p4}. Let $\alpha$ be a color in $\pbar(x)\backslash \{\beta\}$ and $\eta$ be a color in $\pbar(z_1)$. Note that $P_x(\alpha,7,\phiv)$ ends at $t$ by Lemma~\ref{p4link}. Thus we may assume that $\eta\in \pbar(x)$ since otherwise $\eta\in \{1,2,7\}$ and we can interchange $\eta,\alpha$ on $P_{z_1}(\eta,\alpha,\phiv)$. So we assume $\eta\in \pbar(x)$. We then claim that we may further assume that $\eta\in \pbar(x)\backslash \{\phiv(vt)\}$. Otherwise $\eta=\phiv(vt)\in \pbar(x)$.  Interchange $\eta,1$ on $P_{z_1}(\eta,1,\phiv)$ first and then interchange $1,\alpha$ on the $(1,\alpha)$-chain starting at $z_1$. Thus we assume that $\eta\in \pbar(x)\backslash \{\phiv(vt)\}$. Now $P_x(\eta,6,\phiv)$ ends at $z_1$ but not $t$, a contradiction to Lemma~\ref{p4link}.

Now we further assume $\beta \in \{6,7\}$. Without loss of generality assume  $\phiv(vt)=6$. Let $\eta'\in \pbar(t)\backslash \{7\}$.  With a similar argument as before,  we assume $\eta'\in \pbar(x)$. Let $\eta_1$ be the color missing at $z_1$ and $\eta_2$ be the color missing at $z_2$. 

We first claim  $\eta_1=7$. Since otherwise, we have $\eta_1\in \{1,2,3,4,5\}$ and by interchanging $\eta_1,3$ on $P_{z_1}(\eta_1,3,\phiv)$ if necessary, we may assume that $\eta_1\in \pbar(x)$. Then by recoloring $xz_1$ with $\eta_1$, we are back to the case when $\phiv(vt)\in \pbar(x)\cup \pbar(y)$. 

We then claim  $\eta_2=6$. Since otherwise,  $\eta_2\in \{1,2,3,4,5\}$ and by interchanging $\eta_2,3$ on $P_{z_2}(\eta_2,3,\phiv)$ if necessary, we may assume that $\eta_2\in \pbar(x)$. By recoloring $xz_2$ with $\eta_2$ and then recoloring $xz_1$ with $7$, we are back to the case when $\phiv(vt)\in \pbar(x)\cup \pbar(y)$. Thus $\eta_2=6$. Note that the above argument also implies that $P_{z_2}(6,\eta',\phiv)$ ends at $x$, since otherwise by interchanging $6,\eta'$ on this path, we are back to the case when $\eta_2\not=6$.  Now let $\phiv_1=\phiv/P_{t}(\eta',6,\phiv)$, we have $\phiv_1(vt)=\eta'\in \pbar_1(x)$, and thus we are back to the case when $\phiv(vt)\in \pbar(x)\cup \pbar(y)$. This completes the proof of (1). 
\hfill $\Box$

\medskip \noindent
{\bf (2)} 
Since $x$ has no $5^-$-neighbors, by (1) $x$ has two $6$-neighbors and three $7$-neighbors.   Denote by $v_1,v_2,v_3$ the three $7$-vertices   and $z$ the $6$-neighbor of $x$ distinct from $y$. Then $\phiv(xz) \in \{1,2\}$ or  $\phiv(xz) \in \{6,7\}$. 

{\bf Case 2.1} $\phiv(xz)\in \{1,2\}$.

In this case, $x,y,z$ form the vertex set of a multi-fan with respect to $xy$ and $\phiv$. By Lemma~\ref{vf}, we have $\pbar(z)\in \{6,7\}$. Assume without loss of generality that $\phiv(xz)=1$, $\pbar(z)=\{6\}$, $\phiv(xv_1)=2$ and $\phiv(xv_2)=6$. Note that if each of $v_1$ and $v_2$ has at most two $5^-$-neighbors, then we are done. Thus we consider the following two cases.

If  $v_1$ has three $5^-$-neighbors, then there exists $t_1\in N(v_1)\backslash \{x\}$ such that $d(t_1)\le 5$ and $\phiv(v_1t_1)\not=7$. Let $\eta_1$ be a color in $\pbar(t_1)\backslash \{7\}$.  With similar arguments as before we may assume that $\eta_1\in \pbar(x)$ and $\phiv(v_1t_1)\in \pbar(x)$.
 Now $yxv_1t_1$ is a Kierstead path with respect to $xy$ and $\phiv$. But $\eta_1\in \pbar(x)\cap \pbar(t_1)$, a contradiction to Lemma~\ref{p4}.

If $v_2$ has three $5^-$-neighbors, then there exists $t_2\in N(v_2)\backslash \{x\}$ such that $d(t_1)\le 5$ and $\phiv(v_2t_2)\not=7$. Let $\eta_2$ be a color in $\pbar(t_2)\backslash \{7\}$. Similar to the argument before, we may assume that $\eta_2$ and $\phiv(v_2t_2)$ are in $\pbar(x)$. Let $\phiv'=\phiv/P_{t_2}(\eta_2,6,\phiv)$. Then we have $6\in \pbar'(t_2)$. Thus $P_{x}(6,\phiv'(v_2t_2),\phiv')=xv_2t_2$ does not end at $z$,  a contradiction to Lemma~\ref{vf}.  This completes the proof of Case 2.1.

{\bf Case 2.2} $\phiv(xz)\in \{6,7\}.$

In this case, we may assume without loss of generality that $\phiv(xz)=6$, $\phiv(xv_1)=1$ and $\phiv(xv_2)=2$. Note that if each of $v_1$ and $v_2$ has at most two $5^-$-neighbors, then we are done. Thus by the symmetry, assume that $v_1$ has three $5^-$-neighbors. Then there exist two vertices $t,t'\in N(v_1)\backslash \{x\}$ such that $d(t)\le 5$ and $d(t')\le 5$.

{\bf Claim 1} $\{\phiv(v_1t),\phiv(v_1t')\}=\{6,7\}$. 

Otherwise, without loss of generality, assume  $\phiv(v_1t)\in \pbar(x)\cup \pbar(y)$. Then $y,x,v_1,t$ form the vertex set of a Kierstead path with $d(x)<\D$. Thus $\pbar(t)=\{6,7\}$ by Lemma~\ref{p4}. Let $\alpha$ be a color in $\pbar(x)\backslash \{\phiv(v_1t)\}$ and $\eta$ be the color in $\pbar(z)$. Note that $P_x(\alpha,7,\phiv)$ ends at $t$ by Lemma~\ref{p4link}. Thus we may assume that $\eta\in \pbar(x)$ since otherwise $\eta\in \{1,2,7\}$ and we can interchange $\eta,\alpha$ on $P_{z_1}(\eta,\alpha,\phiv)$.  Furthermore, we may assume that $\eta\in \pbar(x)\backslash \{\phiv(v_1t)\}$. Otherwise $\eta=\phiv(v_1t)\in \pbar(x)$, and we can interchange $\eta,1$ on $P_{z}(\eta,1,\phiv)$ first and then interchange $1,\alpha$ on the $(1,\alpha)$-chain starting at $z$.  Now the $(6,\eta)$-chain starting at $x$ ends at $z$ but not $t$,  a contradiction to Lemma~\ref{p4link}. Therefore  $\{\phiv(v_1t),\phiv(v_1t')\}=\{6,7\}$ and without loss of generality, we assume that $\phiv(v_1t)=6$ and $\phiv(v_1t')=7$. This completes the proof of Claim 1.

{\bf Claim 2} $\pbar_1(z)\not=\{7\}$. 

Let $\eta$ be the color missing at $z$. . Otherwise $\eta\in \pbar(x)\cup \pbar(y)$. We may assume that $\eta\in \pbar(x)$ since otherwise we can interchange $\eta$ and $3$ on $P_z(\eta,3,\phiv)$ to get the desired coloring. Now by recoloring $xz$ with $\eta$, we have $\{\phiv(v_1t),\phiv(v_1t')\}\not=\{6,7\}$, a contradiction to Claim 1. Thus $\pbar(z)=\{7\}$.

Now let $\eta'$ be a color in $\pbar(t')\backslash \{6\}$. Similarly as before, we may assume that $\eta'\in \pbar(x)$. If $P_{t'}(\eta',7,\phiv)$ does not end at $x$, let $\phiv_1 = \phiv/ P_{t'}(\eta',7,\phiv)$. Then we have $\{\phiv_1(v_1t),\phiv_2(v_1t')\}\not=\{6,7\}$, a contradiction to Claim 1. If $P_{t'}(\eta',7,\phiv)$ ends at $x$, 
let $\phiv_1 = \phiv/ P_{z}(\eta',7,\phiv)$. Then we have $\pbar_1(z)\not=\{7\}$, a contradiction to Claim 2. This completes the proof of (2). \hfill $\Box$

\medskip \noindent
{\bf  (3)} Since $y$ is the only $6$-neighbor of $x$ and $|\phiv(x)\cap \phiv(y)|=2$, there are two 7-neighbors of $x$, say $v_1,v_2$, such that $\{\phiv(xv_1),\phiv(xv_2)\}\subseteq \pbar(y)$. It is sufficient to show that  each $v_1$ and $v_2$ has at most three $5^-$-neighbors. 

Suppose to the contrary  that $v_1$ has three $5^-$-neighbors other than $x$, say $t_1,t_2,t_3$. Since $|\pbar(x)|=3$, $|\pbar(y)|\ge 2$ and $|\pbar(t_i)|\ge 2$ for each $i=1,2,3$,  by Lemma~\ref{gbroom}, at most one of $\phiv(v_1t_1),\phiv(v_1t_2),\phiv(v_1t_3)$ is in $\pbar(x)\cup \pbar(y)$. Without loss of generality, assume $\phiv(v_1t_1) \in \pbar(x)\cup \pbar(y)$. Then $\{\phiv(v_1t_2),\phiv(v_1t_3)\} = \{6,7\}$.  By Lemma~\ref{p4}, we have $\pbar(t_1)=\phiv(x)\cap \phiv(y)=\{6,7\}$. Thus $\{y,x,v_1,t_1,t_2,t_3\}$ is the vertex set of a $\phiv$-broom. But $\{y,x,v_1,t_1,t_2,t_3\}$ is not elementary, a contradiction to Lemma~\ref{gbroom}. This completes the proof of (3) and thus completes the proof of the lemma.
\end{proof}

\end{document}